\documentclass{article}
\usepackage{amssymb}

\usepackage{graphicx}
\usepackage{amsmath}

\input{tcilatex}

\begin{document}

\title{A new modified Galerkin method for the two-dimensional Navier-Stokes
equations}
\author{Anca-Veronica Ion \\
%EndAName
Faculty of Mathematics and Computer Sciences, \\
University of Pite\c{s}ti, 1,T\^{a}rgu din Vale, Pite\c{s}ti, \\
110040, Arge\c{s}, Romania}
\maketitle

\begin{abstract}
We present a new type of modified Galerkin method. It is a construction with
several (inductively defined) levels, that provides approximate solutions of
increasing accuracy with every new level. These solutions are constructed as
approximations of the so called induced trajectories (notion on which the
definition of a class of approximate inertial manifolds used in the
nonlinear and postprocessed Galerkin methods is based).

Key words: \textit{Navier-Stokes equations, induced trajectories, nonlinear
and postprocessed Galerkin methods, approximate inertial manifolds}

MSC AMS 2000: 35K55, 65M60.
\end{abstract}

\section{Introduction}

We consider the Navier-Stokes equations for a planar flow, with periodic
boundary conditions. The functional framework and the basic hypothesis are
presented in Section 2.

As usual in the Galerkin method, the space of functions is split into a
direct sum of two subspaces: one is the finite dimensional space spanned by
the eigenfunctions corresponding to a finite set, $\Gamma _{m},$ of
eigenvalues of the linear operator $\mathbf{A=-\Delta }$ and the other is
the orthogonal complement of the first. The solution\ $\mathbf{u}$ of the
Navier-Stokes equations will be projected on these spaces: $\mathbf{%
u=Pu+Qu=p+q,}$ where $\mathbf{P}$ is the projector on the finite dimensional
space, and $\mathbf{Q=I-P}$ (Section 3).

In Section 4 we improve the estimates for $\mathbf{q}$ proved in \cite{FMT}.
There, the $\left[ L^{2}\left( \Omega \right) \right] ^{2}$ norm of this
function is found to be less than $K_{0}L^{\frac{1}{2}}\delta ,$ with $%
\delta =\frac{\lambda }{\Lambda }$, where $\lambda $ is the least eigenvalue
of \ $\mathbf{A}$, $\Lambda $ is the least eigenvalue of \ $\mathbf{A}$ not
belonging to $\Gamma _{m},$ and $L$ depends increasingly on the number of
eigenvalues in $\Gamma _{m}.$ The presence of $L$ is not convenient for our
work, since we construct an iterative approximation processus, and if we
would use this estimate, at every step a factor of $L^{1/2}$ would appear in
the evaluation of the error. This \ would \ lead us to bad estimates of the
accuracy of our approximate solutions. We obtain estimates of $\ \mathbf{q}$
independent of $L.$

Our modified Galerkin method is presented in Section 5. The first level of
the method is related to the already classical postprocessed method \cite
{ANT}. This one consists in correcting the Galerkin approximation$\mathbb{\,}
$of the solution, computed at the end of the time integration interval,$\;%
\mathbf{p}_{0}\left( T\right) $ (that is, the solution of (\ref{p0}) at $T$)$%
,$ by adding to it the function$\ \mathbf{q}_{0}\left( T\right) =\Phi _{0}(%
\mathbf{p}_{0}\left( T\right) ).$ $\Phi _{0}(\mathbf{\cdot })$ is the
function whose graph is the approximate inertial manifold (a.i.m.) defined
in \cite{FMT}. The approximation of $\mathbf{u}\left( T\right) $ is taken as 
$\mathbf{p}_{0}\left( T\right) \mathbf{+q}_{0}(T).$ Unlike this
''postprocessed'' Galerkin method, in our method the function $\mathbf{q}%
_{0}\left( t\right) =\Phi _{0}(\mathbf{p}_{0}\left( t\right) )$ is computed
at every moment $t.$ From the numerical point of view this means that it
must be computed at every point of the time grid on $[0,T].$ The approximate
solution at every $t$ is $\mathbf{u}_{0}(t)=\mathbf{p}_{0}\left( t\right) 
\mathbf{+q}_{0}(t).$ The error (in\ $\left[ L_{per}^{2}\left( \Omega \right) %
\right] ^{2}$) of this approximation is of the order of $\delta ^{5/4}.$ We
must remind here the notion of induced trajectory of \cite{T-ind-tr}. There
a family of functions is defined, $\left\{ \mathbf{u}_{j,m}\right\} _{j\geq
0}$, that approximate the exact solution of the Navier-Stokes equations. The
first of these is $\mathbf{u}_{0,m}(t)=\mathbf{p}(t)+\mathbf{q}_{0,m}(t),\;$%
where $\mathbf{q}_{0,m}(t)=\Phi _{0}(\mathbf{p}\left( t\right) ),$ and $%
\mathbf{p}(t)$ is, as above, the $\mathbf{P}$ projection of the exact
solution. So, our function $\mathbf{u}_{0}(t)$ is an approximation of \ $%
\mathbf{u}_{0,m}(t).$

At the second level we look for a new (and better) approximation $\mathbf{p}%
_{1}$ of $\mathbf{p,\;}$by solving an equation closer to the $\mathbf{P}$
projection of the Navier-Stokes equation than the Galerkin equation. That
is, in the equation for $\mathbf{p}_{1},$ in the argument of the nonlinear
term, we approximate $\mathbf{p+q}$ (that appears in the projected N-S
equation (\ref{eqpu})) with $\mathbf{p}_{1}\mathbf{+q}_{0},$ with $\mathbf{q}%
_{0}$ defined above (equation (\ref{def-p1}) in Section 5). The initial
condition for $\mathbf{p}_{1}$ is $\mathbf{Pu}_{0}.\;$The equation (\ref
{def-p1}) is different from those arising in the non-linear Galerkin
methods, since here $\ \mathbf{q}_{0}$ is already determined at the
preceding step, while in the non-linear Galerkin methods \cite{MT}, \cite{DM}
in the equation for $\mathbf{p}$, the argument of the nonlinear term is $%
\mathbf{p+}\Phi _{0}(\mathbf{p})$. Solving equation (\ref{def-p1}) of our
method is not essentially more difficult than solving the Galerkin equation.

Then we compute $\mathbf{q}_{1}\left( t\right) =\widetilde{\Phi }_{1}(%
\mathbf{p}_{1}\left( t\right) ,\mathbf{q}_{0}\left( t\right) )$ where $%
\widetilde{\Phi }_{1}$ is given by (\ref{q1}). The definition of $\mathbf{q}%
_{1}\left( .\right) $ is inspired from that of the function $\mathbf{q}%
_{1,m}\left( .\right) $ of \cite{T-ind-tr} and $\mathbf{q}_{1}\left(
.\right) $ is an approximation of this latter function. There is an obvious
connection between the definition of $\mathbf{q}_{1}\left( t\right) $ and
the second a.i.m. from the family defined in \cite{T-ind-tr} and used in the
non-linear Galerkin methods\ (the definition of this a.i.m. uses that of $%
\mathbf{q}_{1,m}\left( .\right) ).\;$The new approximate solution of the
Navier-Stokes equations we define is $\mathbf{u}_{1}\left( t\right) \mathbf{%
=p}_{1}\left( t\right) \mathbf{+q}_{1}\left( t\right) .$ This is a better
approximation of the exact solution that the preceding one, since the error
is of the order of $\delta ^{7/4}$ (as is shown in Section 6).

We define inductively the next levels of our method. By assuming that we
already found $\mathbf{p}_{j},\;\mathbf{q}_{j}$ for $j\leq k+1$ ($k\geq 0$),
we construct the equation for $\mathbf{p}_{k+2}\mathbf{,}$ taking $\mathbf{p}%
_{k+2}\mathbf{+q}_{k+1}$ in the argument of the nonlinear term (equation (%
\ref{p_(k+2)})). The ''small eddies'' component of the solution will be
approximated by $\mathbf{q}_{k+2}\left( t\right) =\widetilde{\Phi }_{k+2}(%
\mathbf{p}_{k+2}\left( t\right) ,\mathbf{q}_{k+1}\left( t\right) ,\mathbf{q}%
_{k}\left( t\right) ),$ where $\widetilde{\Phi }_{k+2}$ is a function (given
by the right hand side of (\ref{q_k+2})) whose construction was inspired by
that of the function $\mathbf{q}_{k+2,m}\left( t\right) $ of \cite{T-ind-tr}%
. The error of $\mathbf{u}_{k+2}\left( t\right) \mathbf{=p}_{k+2}\left(
t\right) \mathbf{+q}_{k+2}\left( t\right) $ is of the order of $\delta
^{5/4+(k+2)/2}.$ The set $\left\{ \mathbf{u}_{k+2}\left( t\right) \;;\,t\geq
0\right\} $ is an approximation of the induced trajectory $\left\{ \mathbf{u}%
_{k+2,m}\left( t\right) \;;t\geq 0\right\} .$

In Section 6 we prove the estimates of the error of our approximate
solutions. Finally some comments on the advantages and drawbacks of this
method are given in Section 7.

\section{The equations, the functional framework}

The plane flow of an incompressible Newtonian fluid is modelled by the
Navier-Stokes equations:

\begin{eqnarray}
\frac{\partial \mathbf{u}}{\partial t}-\nu \Delta \mathbf{u}+\left( \mathbf{%
u\cdot \nabla }\right) \mathbf{u+}\nabla p &=&\mathbf{f},  \label{eq_u} \\
\text{div}\mathbf{u} &=&0,  \label{inc} \\
\mathbf{u}\left( 0,\mathbf{\cdot }\right) &=&\mathbf{u}_{0}\left( \mathbf{%
\cdot }\right) ,  \label{c-i}
\end{eqnarray}
where $\mathbf{u}=\mathbf{u}\left( t,\mathbf{x}\right) $ is the fluid
velocity, $\mathbf{x\in }\Omega \subset \mathbb{R}^{2},\;\mathbf{u}\left( .,%
\mathbf{x}\right) \mathbf{:}\left[ 0,\infty \right) \rightarrow \mathbb{R}%
^{2}$, $p\left( .,\mathbf{x}\right) :\left[ 0,\infty \right) \rightarrow 
\mathbb{R}$ is the pressure of the fluid, $\nu $ is the kinematic viscosity,
and $\mathbf{f}$ is the volume force. We take here $\Omega =\left(
0,l\right) \times \left( 0,l\right) $ and consider the case of periodic
boundary conditions.

The kinematic viscosity is measured in centistokes ($=mm^{2}/s$). For water
at $20$ Celsius degrees it's value is around 1. In order to have coherent
measure units, we consider the velocity measured in mm/s. We do not focus
here on the mechanics of fluids aspects of the problem, but we focus on the
mathematical construction of the approximate solution. However, we must
remark that the method is appropriate for the study of \ Newtonian fluids
not having very small kinematic viscosity.

We assume that $\mathbf{f}$ is independent of time and is an element of $%
\left[ L_{per}^{2}\left( \Omega \right) \right] ^{2}$. As is usual in the
study of the Navier-Stokes equations with periodic boundary conditions, we
assume that \cite{T-N-S}, \cite{R}

\begin{equation}
\overset{\_}{\mathbf{f}}=\frac{1}{l^{2}}\int_{\Omega }\mathbf{f}\left( 
\mathbf{x}\right) d\mathbf{x=0,}  \label{f}
\end{equation}
and that the pressure is a periodic function on $\Omega $. For simplicity we
will assume also that the average of the velocity over the periodicity cell
is zero.

The velocity $\mathbf{u}$\textbf{\ }is thus looked for in the space $\ $%
\newline
$\mathcal{H}=\left\{ \mathbf{v;\;v\in }\left[ L_{per}^{2}\left( \Omega
\right) \right] ^{2},\;\text{div}\,\mathbf{v=0,\;}\overset{\_}{\mathbf{u}}%
=0\right\} $. The scalar product in $\ \mathcal{H}$ is $\ \left( \mathbf{u},%
\mathbf{v}\right) =\int_{\Omega }\left( u_{1}v_{1}+u_{2}v_{2}\right) dx,$
(where $\mathbf{u}=\left( u_{1},u_{2}\right) ,\;\mathbf{v}=\left(
v_{1},v_{2}\right) $). $\;$The induced norm is denoted by $\left| \quad
\right| $.

We also need the space $\mathcal{V}=\left\{ \mathbf{u}\in \left[
H_{per}^{1}\left( \Omega \right) \right] ^{2},\;div\,\mathbf{u=0,}\overset{\_%
}{\mathbf{u}}=0\right\} ,$ with the scalar product $\left( \left( \mathbf{u},%
\mathbf{v}\right) \right) =\sum_{i,j=1}^{2}\left( \frac{\partial u_{i}}{%
\partial x_{j}},\frac{\partial v_{i}}{\partial x_{j}}\right) ,$ and the
induced norms, denoted by $\left\| \;\right\| .$ We denote $\mathbf{%
A=-\Delta }$ and observe that $\mathbf{A}$ is defined on $D(\mathbf{A})=%
\mathcal{H\cap }H^{2}\left( \Omega \right) .$

We shall focus on finding approximations for the function $\mathbf{u}$.

The classical variational formulation of the Navier-Stokes equations \cite
{T-N-S} leads to the abstract equation 
\begin{eqnarray}
\frac{d\mathbf{u}}{dt}-\nu \Delta \mathbf{u}+\left( \mathbf{u\cdot \nabla }%
\right) \mathbf{u} &=&\mathbf{f}\;\;\;\text{in\ }\mathcal{V\,}^{\prime },
\label{u-tilda} \\
\mathbf{u}\left( 0\right) &=&\mathbf{u}_{0},\;\;\mathbf{u}_{0}\in \mathcal{H}%
.  \label{ci-u-tilda}
\end{eqnarray}

The notations $\mathbf{B(u,v)=(u\cdot \nabla )v}$, $\mathbf{%
B(u)=B(u,u),\;b(u,v,w)=}\left( \mathbf{B(u,v),w}\right) \;$

\noindent will be used below.

For the bilinear application $\mathbf{B}(\mathbf{u,v})$ the following
inequalities

\begin{eqnarray}
\left| \mathbf{B}\left( \mathbf{u,v}\right) \right| &\leq &c_{1}\left| 
\mathbf{u}\right| ^{\frac{1}{2}}\left| \mathbf{\Delta u}\right| ^{\frac{1}{2}%
}\left\| \mathbf{v}\right\| ,\text{ \ \ \ \ \ \ \ \ \ \ \ \ \ \ \ \ \ }%
\left( \forall \right) \mathbf{\ u\in }D(\mathbf{A}),\;\mathbf{v\in }%
\mathcal{V},  \label{B(u,v)1} \\
\ \left| \mathbf{B}\left( \mathbf{u,v}\right) \right| &\leq &c_{2}\left\| 
\mathbf{u}\right\| \left\| \mathbf{v}\right\| \left[ 1+\ln \left( \frac{%
\left| \mathbf{\Delta u}\right| ^{2}}{\lambda _{1}\left\| \mathbf{u}\right\|
^{2}}\right) \right] ^{\frac{1}{2}},\;\left( \forall \right) \mathbf{\ u\in }%
D(\mathbf{A}),\;\mathbf{v\in }\mathcal{V}.  \label{B(u,v)2}
\end{eqnarray}
hold \cite{AG}, \cite{T-N-S}, \cite{T-ind-tr}. We remind the following
properties of the trilinear form $\mathbf{b}(\mathbf{u},\mathbf{v},\mathbf{w}%
)$ (valid for periodic boundary conditions \cite{R}):$\ \ \ \ \ \ \ \ \ $%
\begin{eqnarray}
\ \mathbf{b}(\mathbf{u},\mathbf{v},\mathbf{w}) &=&-\mathbf{b}(\mathbf{u},%
\mathbf{w},\mathbf{v}),  \label{b0} \\
\mathbf{b}(\mathbf{u},\mathbf{v},\mathbf{v}) &=&0,  \label{b1}
\end{eqnarray}

\noindent as well as the following inequalities \cite{R}

\begin{eqnarray}
\left| \mathbf{b}(\mathbf{u},\mathbf{v},\mathbf{w})\right| &\leq
&c_{3}\left| \mathbf{u}\right| ^{\frac{1}{2}}\left\| \mathbf{u}\right\| ^{%
\frac{1}{2}}\left\| \mathbf{v}\right\| \left| \mathbf{w}\right| ^{\frac{1}{2}%
}\left\| \mathbf{w}\right\| ^{\frac{1}{2}},\;\left( \forall \right) \mathbf{%
u,v,w\in }\mathcal{V},  \label{b(u,v,w)1} \\
\ \left| \mathbf{b}(\mathbf{u},\mathbf{v},\mathbf{w})\right| &\leq
&c_{4}\left| \mathbf{u}\right| ^{\frac{1}{2}}\left\| \mathbf{u}\right\| ^{%
\frac{1}{2}}\left\| \mathbf{v}\right\| ^{\frac{1}{2}}\left| \Delta \mathbf{v}%
\right| ^{\frac{1}{2}}\left| \mathbf{w}\right| ,\;\left( \forall \right) 
\mathbf{u\in }\mathcal{V},\;\mathbf{v\in D}\left( \mathbf{A}\right) ,\mathbf{%
w\in }\mathcal{H}.\qquad  \label{b(u,v,w)2}
\end{eqnarray}

For the problem (\ref{u-tilda})\textit{, }(\ref{ci-u-tilda}) we have the
classical existence and uniqueness results for the equations Navier-Stokes
in $\mathbb{R}^{2}$, with periodic boundary conditions.

\textbf{Theorem 1 }\cite{T-N-S}. \textit{a)} \textit{If \ }$\mathbf{u}%
_{0}\in \mathcal{H},\;\mathbf{f}\in \mathcal{H},$\textit{\ then the problem }
(\ref{u-tilda})\textit{, }(\ref{ci-u-tilda}) \textit{has an unique solution}%
\textbf{\ }\textit{\ }$\mathbf{u}\in C^{0}\left( [0,T];\mathcal{H}\right)
\cap L^{2}\left( 0,T;\mathcal{V}\right) .\;b)$ \textit{If, \ in addition to
the hypotheses in a), } $\mathbf{u}_{0}\in \mathcal{V},\;$\textit{then\ }%
\textbf{\ }\textit{\ }$\mathbf{u}\in C^{0}\left( [0,T];\mathcal{V}\right)
\cap L^{2}\left( 0,T;D(\mathbf{A})\right) .$\textit{\ The solution is, in
this latter case, analytic in time on the positive real axis.}

The semi-dynamical system $\left\{ S\left( t\right) \right\} _{t\geq 0}$
generated by problem (\ref{u-tilda}) is dissipative \cite{T}, \cite{R}.%
\textbf{\ }More precisely, there is a $\rho _{0}>0$ $\ $ such that for every 
$R>0,$ there is a $t_{0}(R)>0$ with the property that for every $\mathbf{u}%
_{0}\in \mathcal{H}\;$with\textbf{\ } $\left| \mathbf{u}_{0}\right| \leq R$,
we have $\left| S\left( t\right) \mathbf{u}_{0}\right| \leq \rho _{0}$ for $%
t>t_{0}(R)$ $.$ In addition, there are absorbing balls in$\;\mathcal{V\;}$%
and $\mathbf{D}\left( \mathbf{A}\right) $ for $\left\{ S\left( t\right)
\right\} _{t\geq 0}$, i.e. there are $\rho _{1}>0,$ $\rho _{2}>0$ and $%
t_{1}(R),\;t_{2}(R)$ with $t_{2}(R)\geq t_{1}(R)\geq t_{0}(R)$ such that for
every $R>0,\;$ $\left| \mathbf{u}_{0}\right| \leq R$ implies $\left\|
S\left( t\right) \mathbf{u}_{0}\right\| \leq \rho _{1}$ for $t>t_{1}(R)$ \
and\ $\left| \mathbf{A}S\left( t\right) \mathbf{u}_{0}\right| \leq \rho _{2}$
for $t>t_{2}(R).$

\section{The decomposition of the space, the projections of the equations}

The eigenvalues of $\mathbf{A}$ are $\lambda _{j_{1},j_{2}}=\frac{4\pi ^{2}}{%
l^{2}}\left( j_{1}^{2}+j_{2}^{2}\right) ,\;\;\left( j_{1},\;j_{2}\right) \in 
\mathbb{N}^{2}\backslash \left\{ \left( 0,0\right) \right\} ,$ and the
corresponding eigenfunctions are 
\begin{eqnarray*}
\mathbf{w}_{j_{1},j_{2}}^{s\pm } &=&\frac{\sqrt{2}}{l}\frac{\left( j_{2},\mp
j_{1}\right) }{\left| \mathbf{j}\right| }\sin \left( 2\pi \frac{%
j_{1}x_{1}\pm j_{2}x_{2}}{l}\right) , \\
\mathbf{w}_{j_{1},j_{2}}^{c\pm } &=&\;\frac{\sqrt{2}}{l}\frac{\left(
j_{2},\mp j_{1}\right) }{\left| \mathbf{j}\right| }\cos \left( 2\pi \frac{%
j_{1}x_{1}\pm j_{2}x_{2}}{l}\right) ,
\end{eqnarray*}
where $\left| \mathbf{j}\right| =\left( j_{1}^{2}+j_{2}^{2}\right) ^{\frac{1%
}{2}}$\cite{T-ind-tr}$.$ These eigenfunctions form a total system for $%
\mathcal{H}$. In order to be able to write easily sums involving the four
eigenfunctions above, we denote them as follows 
\begin{equation*}
\mathbf{w}_{j_{1},j_{2}}^{s+}=\mathbf{w}_{j_{1},j_{2}}^{1},\;\;\mathbf{w}%
_{j_{1},j_{2}}^{s-}=\mathbf{w}_{j_{1},j_{2}}^{2},\;\;\mathbf{w}%
_{j_{1},j_{2}}^{c+}=\mathbf{w}_{j_{1},j_{2}}^{3},\;\;\mathbf{w}%
_{j_{1},j_{2}}^{c-}=\mathbf{w}_{j_{1},j_{2}}^{4}.
\end{equation*}

For a fixed $m\in \mathbb{N}$ we consider the set $\Gamma _{m}\;$of
eigenvalues $\lambda _{j_{1},j_{2}}$ having $\ 0\leq j_{1},j_{2}\leq m$. We
define$\;$ 
\begin{eqnarray*}
\lambda &:&=\lambda _{1,0}=\lambda _{0,1}=\frac{4\pi ^{2}}{l^{2}}, \\
\Lambda &:&=\lambda _{m+1,0}=\lambda _{0,m+1}=\frac{4\pi ^{2}}{l^{2}}\left(
m+1\right) ^{2}, \\
\delta &=&\delta \left( m\right) :=\frac{\lambda }{\Lambda }=\frac{1}{%
(m+1)^{2}}.
\end{eqnarray*}

$\Lambda $ is the least eigenvalue not belonging to $\Gamma _{m}$. The
eigenfunctions corresponding to the eigenvalues of $\ \Gamma _{m}$ span a
finite-dimensional subspace of $\mathcal{H}$. We denote by $\mathbf{P}$ the
orthogonal projection operator on this subspace and by $\mathbf{Q}$ the
orthogonal projection operator on the complementary subspace. We write for
the solution $\mathbf{u}$ of (\ref{u-tilda}),

\begin{center}
$\mathbf{p=Pu},\;\mathbf{q=Qu}.$
\end{center}

By projecting equation (\ref{u-tilda}) on the above constructed spaces, we
obtain 
\begin{eqnarray}
\frac{d\mathbf{p}}{dt}-\nu \Delta \mathbf{p}+\mathbf{PB}(\mathbf{p+q}) &=&%
\mathbf{Pf},  \label{eqpu} \\
\frac{d\mathbf{q}}{dt}-\nu \Delta \mathbf{q}+\mathbf{QB}(\mathbf{p+q}) &=&%
\mathbf{Qf}.  \label{eqqu}
\end{eqnarray}

\section{New estimates for the ''small'' component of the solution}

In \cite{FMT} is proved that for every $R>0,$ there is a moment $t_{3}\left(
R\right) \geq t_{2}(R)$ such that for every $\left| \mathbf{u}_{0}\right|
\leq R,$ 
\begin{eqnarray}
\left| \mathbf{q}\left( t\right) \right| &\leq &K_{0}L^{\frac{1}{2}}\delta
,\;\;\;\;\left\| \mathbf{q}\left( t\right) \right\| \leq K_{1}L^{\frac{1}{2}%
}\delta ^{\frac{1}{2}},  \label{est-q} \\
\left| \mathbf{q}^{\prime }\left( t\right) \right| &\leq &K_{0}^{\prime }L^{%
\frac{1}{2}}\delta ,\;\;\;\;\;\left| \Delta \mathbf{q}\left( t\right)
\right| \leq K_{2}L^{\frac{1}{2}},\;\;\;\;\;\ t\geq t_{3}\left( R\right) , 
\notag
\end{eqnarray}
where $K_{0},$ $K_{0}^{\prime },\;K_{1},$ $K_{2}$ depend of $\nu ,\;\left|
f\right| ,\;\lambda \;$and,$\;$for the way we chose the projection
subspaces, $L=L(m)=1+\ln 2m^{2}$ (see also \cite{T-ind-tr}). The constant $L$
comes from the use of inequality (\ref{B(u,v)2}) in the course of the proof
of (\ref{est-q}). More specific 
\begin{eqnarray*}
L &=&\underset{\mathbf{p\in P}\mathcal{H}}{\sup }\left( 1+\ln \frac{\left|
\Delta \mathbf{p}\right| ^{2}}{\lambda _{1}\left\| \mathbf{p}\right\| ^{2}}%
\right) =\underset{\lambda \in \Gamma _{m}}{\max }\left( 1+\ln \frac{\lambda 
}{\lambda _{1}}\right) = \\
&=&1+\ln 2m^{2}.
\end{eqnarray*}

In the sequel we shall improve the above estimates, trying to eliminate $L$
(which tends to infinity with $m)$ from the constants. The idea is that of
refining the contribution of the term $\mathbf{QB}(\mathbf{p})$ resulting
from $\mathbf{QB}(\mathbf{p+q})$ in (\ref{eqqu}). We start from the
trigonometric relation 
\begin{eqnarray*}
\sin \left( 2\pi \frac{j_{1}x_{1}\pm j_{2}x_{2}}{l}\right) \sin \left( 2\pi 
\frac{k_{1}x_{1}\pm k_{2}x_{2}}{l}\right) &=&\frac{1}{2}\left[ \cos 2\pi 
\frac{\left( j_{1}-k_{1}\right) x_{1}\pm \left( j_{2}-k_{2}\right) x_{2}}{l}%
\right. - \\
&&-\left. \cos 2\pi \frac{\left( j_{1}+k_{1}\right) x_{1}\pm \left(
j_{2}+k_{2}\right) x_{2}}{l}\right] ,
\end{eqnarray*}
and the similar ones for all other combinations of sine and cosine that
might appear in the scalar product of two eigenfunctions. Since $\mathbf{p=}%
\sum\limits_{0\leq j_{1},j_{2}\leq m}\sum\limits_{i=1}^{4}p_{j_{1},j_{2}}^{i}%
\mathbf{w}_{j_{1},j_{2}}^{i},$ from $\left( \mathbf{p\nabla }\right) \mathbf{%
p=}\left( \sum\limits_{0\leq j_{1},j_{2}\leq
m}\sum\limits_{i=1}^{4}p_{j_{1},j_{2}}^{i}\mathbf{w}_{j_{1},j_{2}}^{i}\nabla
\right) \left( \sum\limits_{0\leq k_{1},k_{2}\leq
m}\sum\limits_{l=1}^{4}p_{k_{1},k_{2}}^{l}\mathbf{w}_{k_{1},k_{2}}^{l}%
\right) ,$ only those products of terms that have $j_{1}+k_{1}\geq m+1$ or $%
j_{2}+k_{2}\geq m+1\;$will belong to $\mathbf{Q}\mathcal{H}.$

We consider from this point on that $m$ is even, and we set $m=2n.$

If for $\mathbf{w}_{j_{1},j_{2}}^{i}$ and $\mathbf{w}_{k_{1},k_{2}}^{l}\;$we
have $j_{1},j_{2}\leq n$ and $k_{1},k_{2}\leq n,$ then \newline
$\left( \mathbf{w}_{j_{1},j_{2}}^{i}\nabla \right) $ $\mathbf{w}%
_{k_{1},k_{2}}^{l}$ belongs to $\mathbf{P}\mathcal{H}$. We are led to the
idea of considering the subspace of $\mathcal{H}$ spanned by all the
eigenfunctions $\mathbf{w}_{j_{1},j_{2}}^{i}$ with $0\leq $ $j_{1},$ $j_{2}$ 
$\leq n,$ $1\leq i\leq 4$. We denote by $\mathbf{P}_{\mathbf{p}}$ the
projection operator on this space and set $\mathbf{P}_{\mathbf{q}}=\mathbf{%
P-P}_{\mathbf{p}},$ $\mathbf{p}_{p}=\mathbf{P}_{\mathbf{p}}\mathbf{p}$ and $%
\mathbf{p}_{q}=\mathbf{P}_{\mathbf{q}}\mathbf{p}$. Obviously 
\begin{equation*}
\mathbf{Q}\left( \mathbf{p}_{p}\mathbf{\nabla }\right) \mathbf{p}_{p}\mathbf{%
=0.}
\end{equation*}

On another hand, we see that $\mathbf{p}_{q}$ is a truncation of $(\mathbf{%
I-P}_{\mathbf{p}})\mathbf{u}$, hence \newline
$\left| \mathbf{p}_{q}\right| \leq \left| (\mathbf{I-P}_{\mathbf{p}})\mathbf{%
u}\right| ,\;\left\| \mathbf{p}_{q}\right\| \leq \left\| (\mathbf{I-P}_{%
\mathbf{p}})\mathbf{u}\right\| .$ Then, by setting $\delta _{1}=\delta
\left( n\right) =\frac{1}{(n+1)^{2}},$ $L_{1}=L(n)=1+\ln 2n^{2}$, the
estimates (\ref{est-q}) imply 
\begin{eqnarray}
\left| \mathbf{p}_{q}\right| &\leq &K_{0}L_{1}^{\frac{1}{2}}\delta
_{1},\;\;\;\;\left\| \mathbf{p}_{q}\right\| \leq K_{1}L_{1}^{\frac{1}{2}%
}\delta _{1}^{\frac{1}{2}},\;  \label{p_q} \\
\left| \mathbf{p}_{q}^{\prime }\right| &\leq &K_{0}^{\prime }L_{1}^{\frac{1}{%
2}}\delta _{1},\;\;\;\;\;\left| \Delta \mathbf{p}_{q}\right| \leq
K_{2}L_{1}^{\frac{1}{2}}.  \label{p_q1}
\end{eqnarray}
We use these inequalities in order to refine the estimates (\ref{est-q}). In
the rest of the paper we shall assume that for a fixed $R\geq 0,$ the
function $\mathbf{u}_{0}$ is such that $\left| \mathbf{u}_{0}\right| \leq R$%
. We state and prove

\textbf{Theorem 1. }\textit{There are some constants }$\widetilde{C}_{0},\;%
\widetilde{C}_{1},\;\widetilde{C}\,_{0}^{\prime },\;\widetilde{C}_{2},$ 
\textit{depending only} $\ $\textit{on }$\nu ,$ $\lambda ,\;\left| \mathbf{Qf%
}\right| $ \textit{such that, for} $t$\textit{\ large enough,} \textit{the
inequalities} 
\begin{eqnarray}
\left| \mathbf{q(t)}\right| &\leq &\widetilde{C}_{0}\delta ,  \label{ref q}
\\
\left\| \mathbf{q(t)}\right\| &\leq &\widetilde{C}_{1}\delta ^{\frac{1}{2}},
\label{ref norm1q} \\
\left| \mathbf{q}^{\prime }\mathbf{(t)}\right| &\leq &\widetilde{C}%
\,_{0}^{\prime }\delta ,  \label{ref q prim} \\
\left| \Delta \mathbf{q(t)}\right| &\leq &\widetilde{C}_{2},
\label{ref delta q}
\end{eqnarray}
\textit{hold.}

\textbf{Proof. }We have, with the notation settled before the Proposition, 
\begin{eqnarray}
\mathbf{QB}(\mathbf{p}) &=&\mathbf{QB}(\mathbf{p}_{p}\mathbf{+p}_{q})=%
\mathbf{QB}(\mathbf{p}_{p})+\mathbf{QB}(\mathbf{p}_{p}\mathbf{,p}_{q})+%
\mathbf{QB}(\mathbf{p}_{q}\mathbf{,p}_{p})+\mathbf{QB}(\mathbf{p}_{q}) 
\notag \\
&=&\mathbf{QB}(\mathbf{p}_{p}\mathbf{,p}_{q})+\mathbf{QB}(\mathbf{p}_{q}%
\mathbf{,p}_{p})+\mathbf{QB}(\mathbf{p}_{q}),  \label{desc B}
\end{eqnarray}
since $\mathbf{QB}(\mathbf{p}_{p})=0.$ Now, as is usual, we take the scalar
product of (\ref{eqqu}) with $\mathbf{q,}$ and by using (\ref{b1}) and (\ref
{desc B}), we obtain 
\begin{eqnarray}
\frac{1}{2}\frac{d\left| \mathbf{q}\right| ^{2}}{dt}+\nu \left\| \mathbf{q}%
\right\| ^{2} &\leq &\left| \left( \mathbf{B}(\mathbf{p}_{p}\mathbf{,p}_{q}),%
\mathbf{q}\right) \right| +\left| \left( \mathbf{B}(\mathbf{p}_{q}\mathbf{,p}%
_{p}),\mathbf{q}\right) \right| +  \notag \\
&&+\left| \left( \mathbf{B}(\mathbf{p}_{q}),\mathbf{q}\right) \right|
+\left| \left( \mathbf{B}(\mathbf{q,p}),\mathbf{q}\right) \right| +\left|
\left( \mathbf{Qf,q}\right) \right| .  \label{ineg-q}
\end{eqnarray}
For the first term of the right-hand side, the following estimates (obtained
by using (\ref{B(u,v)1}) and (\ref{p_q})) hold

\begin{eqnarray*}
\left| \left( \mathbf{B}(\mathbf{p}_{p}\mathbf{,p}_{q}),\mathbf{q}\right)
\right| &\leq &c_{1}\left| \mathbf{p}_{p}\right| ^{1/2}\left| \Delta \mathbf{%
p}_{p}\right| ^{1/2}\left\| \mathbf{p}_{q}\right\| \left| \mathbf{q}\right|
\leq c_{1}\rho _{0}^{1/2}\rho _{2}^{1/2}K_{1}L_{1}^{\frac{1}{2}}\delta _{1}^{%
\frac{1}{2}}\frac{1}{\Lambda ^{\frac{1}{2}}}\text{ }\left\| \mathbf{q}%
\right\| \\
&\leq &c_{1}^{2}\rho _{0}\rho _{2}K_{1}^{2}L_{1}\delta _{1}\frac{2}{\nu
\Lambda }+\frac{\nu }{8}\left\| \mathbf{q}\right\| ^{2}.
\end{eqnarray*}

For the second term we obtain the inequalities 
\begin{eqnarray*}
\left| \left( \mathbf{B}(\mathbf{p}_{q}\mathbf{,p}_{p}),\mathbf{q}\right)
\right| &\leq &c_{1}\left| \mathbf{p}_{q}\right| ^{1/2}\left| \Delta \mathbf{%
p}_{q}\right| ^{1/2}\left\| \mathbf{p}_{p}\right\| \left| \mathbf{q}\right|
\leq c_{1}K_{1}^{1/2}L_{1}^{\frac{1}{4}}\delta _{1}^{\frac{1}{2}}\rho
_{2}^{1/2}\rho _{1}\frac{1}{\Lambda ^{\frac{1}{2}}}\text{ }\left\| \mathbf{q}%
\right\| \\
&\leq &c_{1}^{2}K_{1}\rho _{1}^{2}\rho _{2}L_{1}^{1/2}\delta _{1}\text{ }%
\frac{2}{\nu \Lambda }+\frac{\nu }{8}\left\| \mathbf{q}\right\| ^{2}.
\end{eqnarray*}

\bigskip For the third term we have 
\begin{eqnarray*}
\left| \left( \mathbf{B}(\mathbf{p}_{q}),\mathbf{q}\right) \right| &\leq
&c_{2}L^{\frac{1}{2}}\left\| \mathbf{p}_{q}\right\| ^{2}\left| \mathbf{q}%
\right| \leq c_{2}L^{\frac{1}{2}}K_{1}^{2}L_{1}\delta _{1}\frac{1}{\Lambda ^{%
\frac{1}{2}}}\text{ }\left\| \mathbf{q}\right\| \\
&\leq &c_{2}^{2}K_{1}^{4}LL_{1}^{2}\delta _{1}^{2}\frac{2}{\nu \Lambda }+%
\frac{\nu }{8}\left\| \mathbf{q}\right\| ^{2},
\end{eqnarray*}
and for the fourth, by using (\ref{b(u,v,w)2}) 
\begin{eqnarray*}
\left| \left( \mathbf{B}(\mathbf{q,p}),\mathbf{q}\right) \right| &\leq
&c_{4}\left| \mathbf{q}\right| ^{\frac{1}{2}}\left\| \mathbf{q}\right\| ^{%
\frac{1}{2}}\left\| \mathbf{p}\right\| ^{\frac{1}{2}}\left| \Delta \mathbf{p}%
\right| ^{\frac{1}{2}}\left| \mathbf{q}\right| \leq \\
&\leq &c_{4}K_{0}^{1/2}L^{\frac{1}{4}}\delta ^{\frac{1}{2}}K_{1}^{1/2}L^{%
\frac{1}{4}}\delta ^{\frac{1}{4}}\rho _{1}^{\frac{1}{2}}\rho _{2}^{\frac{1}{2%
}}\frac{1}{\Lambda ^{\frac{1}{2}}}\text{ }\left\| \mathbf{q}\right\| \\
&\leq &\frac{2}{\nu \Lambda }c_{4}^{2}K_{0}K_{1}\rho _{1}\rho _{2}L\delta ^{%
\frac{3}{2}}+\frac{\nu }{8}\left\| \mathbf{q}\right\| ^{2}.
\end{eqnarray*}

At last 
\begin{equation*}
\left| \left( \mathbf{Qf,q}\right) \right| \leq \left| \mathbf{Qf}\right|
\left| \mathbf{q}\right| \leq \frac{2\left| \mathbf{Qf}\right| ^{2}}{\nu
\Lambda }+\frac{\nu }{8}\left\| \mathbf{q}\right\| ^{2}.
\end{equation*}

The above inequalities and (\ref{ineg-q}) lead us to 
\begin{equation*}
\frac{1}{2}\frac{d}{dt}\left| \mathbf{q}\right| ^{2}+3\frac{\nu \Lambda }{8}%
\left| \mathbf{q}\right| ^{2}\leq C_{0}^{2}\delta ,
\end{equation*}
with\ 

\begin{eqnarray*}
C_{0}^{2} &=&\frac{2}{\nu \lambda }\left[ c_{1}^{2}\rho _{0}\rho
_{2}K_{1}^{2}L_{1}\delta _{1}+c_{1}^{2}K_{1}\rho _{1}^{2}\rho
_{2}L_{1}^{1/2}\delta _{1}+c_{2}^{2}K_{1}^{4}LL_{1}^{2}\delta
_{1}^{2}+\right. \\
&&\left. +c_{4}^{2}K_{0}K_{1}\rho _{1}\rho _{2}L\delta ^{\frac{3}{2}}+\left| 
\mathbf{Qf}\right| ^{2}\right] .
\end{eqnarray*}

It follows, with the usual Gronwall Lemma, 
\begin{equation*}
\left| \mathbf{q}\left( t\right) \right| ^{2}\leq \left| \mathbf{q}\left(
0\right) \right| ^{2}e^{-\frac{3}{4}\nu \Lambda t}+\frac{8C_{0}^{2}}{3\nu
\lambda }\delta ^{2},
\end{equation*}
hence, for $t_{4}\left( R\right) \geq t_{3}(R)$, taken as to have $\left| 
\mathbf{q}\left( 0\right) \right| ^{2}e^{-\frac{3}{4}\nu \Lambda t}\leq 
\frac{8C_{0}^{2}}{3\nu \lambda }\delta ^{2}$ for \newline
$t\geq t_{4}\left( R\right) ,$ we obtain (\ref{ref q}), with $\widetilde{C}%
_{0}=\frac{4C_{0}}{\sqrt{3\nu \lambda }}.$

The functions of $n:\;L_{1}\delta _{1}=L\left( n\right) \delta \left(
n\right) =\left( 1+\ln 2n^{2}\right) /\left( n+1\right) ^{2},$ \newline
$L_{1}^{1/2}\delta _{1}=\sqrt{L\left( n\right) }\delta \left( n\right) =%
\sqrt{1+\ln 2n^{2}}/\left( n+1\right) ^{2}$, \newline
$LL_{1}^{2}\delta _{1}^{2}=L\left( 2n\right) L\left( n\right) ^{2}\delta
\left( n\right) ^{2}=\left( 1+\ln 8n^{2}\right) \left( 1+\ln 2n^{2}\right)
^{2}/\left( n+1\right) ^{4}$ and \newline
$L\delta ^{\frac{3}{2}}=\left( 1+\ln 8n^{2}\right) /\left( 2n+1\right)
^{3},\;$that appear in the structure of $C_{0}^{2},$ have at $n=2$ values
less than 1 and are decreasing when $n$ increases (for $n\geq 2).$ Then, for 
$n\geq 2$ 
\begin{equation*}
C_{0}^{2}\leq \frac{2}{\nu \lambda }\left( c_{1}^{2}\rho _{0}\rho
_{2}K_{1}^{2}+c_{1}^{2}K_{1}\rho _{1}^{2}\rho
_{2}+c_{2}^{2}K_{1}^{4}+c_{4}^{2}K_{0}K_{1}\rho _{1}\rho _{2}+\left| \mathbf{%
Qf}\right| ^{2}\right)
\end{equation*}
and the right hand side depends only on $\nu ,\;\lambda ,\;\left| \mathbf{Qf}%
\right| .$ More than that, since all the functions defined above tend to
zero when $n\rightarrow \infty $, we can choose $n$ large enough so that $%
\frac{2\left| \mathbf{Qf}\right| ^{2}}{\nu \lambda }$ becomes the dominant
term in $C_{0}^{2}$.

For that $n$, $\widetilde{C}_{0}$ will be of the order of $\frac{\left| 
\mathbf{Qf}\right| }{\nu \lambda }$ . However, the structure of $%
K_{0},\;K_{1},\;\rho _{0},\;\rho _{1},\;\rho _{2}$ show that if $\ \nu \;$is
very small, then $n$ with the above property must be very large.

\bigskip

Now, we aim to estimate $\left\| \mathbf{q}\right\| .$ By multiplying
equation (\ref{eqqu}) by $\Delta \mathbf{q}$, and by using (\ref{desc B}),
we obtain 
\begin{eqnarray*}
\frac{1}{2}\frac{d\left\| \mathbf{q}\right\| ^{2}}{dt}+\nu \left| \Delta 
\mathbf{q}\right| ^{2} &\leq &\left| \left( \mathbf{B}(\mathbf{p}_{p}\mathbf{%
,p}_{q}+\mathbf{q}),\Delta \mathbf{q}\right) \right| +\left| \left( \mathbf{B%
}(\mathbf{p}_{q}\mathbf{+q,p}_{p}),\Delta \mathbf{q}\right) \right| + \\
&&+\left| \left( \mathbf{B}(\mathbf{p}_{q}\mathbf{+q,p}_{q}\mathbf{+q}%
),\Delta \mathbf{q}\right) \right| + \\
&&+\left| \left( \mathbf{Qf,\Delta q}\right) \right| .
\end{eqnarray*}

For the first term in the right hand side we have (\ref{B(u,v)1}) 
\begin{eqnarray*}
\left| \left( \mathbf{B}(\mathbf{p}_{p}\mathbf{,p}_{q}+\mathbf{q}),\Delta 
\mathbf{q}\right) \right| &\leq &c_{1}\left| \mathbf{p}_{p}\right| ^{\frac{1%
}{2}}\left| \mathbf{\Delta p}_{p}\right| ^{\frac{1}{2}}\left\| \mathbf{p}%
_{q}+\mathbf{q}\right\| \left| \Delta \mathbf{q}\right| \\
&\leq &c_{1}\rho _{0}^{1/2}\rho _{2}^{1/2}L_{1}^{\frac{1}{2}}K_{1}\delta
_{1}^{\frac{1}{2}}\left| \Delta \mathbf{q}\right| \\
&\leq &c_{1}^{2}\rho _{0}\rho _{2}K_{1}^{2}L_{1}\delta _{1}\frac{2}{\nu }+%
\frac{\nu }{8}\left| \Delta \mathbf{q}\right| ^{2},
\end{eqnarray*}
for the second, with (\ref{b(u,v,w)2}), 
\begin{eqnarray*}
\left| \left( \mathbf{B}(\mathbf{p}_{q}\mathbf{+q,p}_{p}),\Delta \mathbf{q}%
\right) \right| &\leq &c_{4}\left| \mathbf{p}_{q}\mathbf{+q}\right| ^{\frac{1%
}{2}}\left\| \mathbf{p}_{q}\mathbf{+q}\right\| ^{\frac{1}{2}}\left\| \mathbf{%
p}_{p}\right\| ^{\frac{1}{2}}\left| \Delta \mathbf{p}_{p}\right| ^{\frac{1}{2%
}}\left| \Delta \mathbf{q}\right| \\
&\leq &c_{4}\widetilde{C}_{0}^{\frac{1}{2}}\delta _{1}^{\frac{1}{2}}L_{1}^{%
\frac{1}{4}}K_{1}^{\frac{1}{2}}\delta _{1}^{\frac{1}{4}}\rho _{1}^{\frac{1}{2%
}}\rho _{2}^{\frac{1}{2}}\left| \Delta \mathbf{q}\right| \\
&\leq &c_{4}^{2}\widetilde{C}_{0}\rho _{1}\rho _{2}K_{1}L_{1}^{\frac{1}{2}%
}\delta _{1}^{\frac{3}{2}}\frac{2}{\nu }+\frac{\nu }{8}\left| \Delta \mathbf{%
q}\right| ^{2}
\end{eqnarray*}
for the third, also with (\ref{b(u,v,w)2}), 
\begin{eqnarray*}
\left| \left( \mathbf{B}(\mathbf{p}_{q}\mathbf{+q,p}_{q}\mathbf{+q}),\Delta 
\mathbf{q}\right) \right| &\leq &c_{4}\left| \mathbf{p}_{q}\mathbf{+q}%
\right| ^{\frac{1}{2}}\left\| \mathbf{p}_{q}\mathbf{+q}\right\| \left|
\Delta \left( \mathbf{p}_{q}\mathbf{+q}\right) \right| ^{\frac{1}{2}}\left|
\Delta \mathbf{q}\right| \\
&\leq &c_{4}\widetilde{C}_{0}^{\frac{1}{2}}\delta _{1}^{\frac{1}{2}}L_{1}^{%
\frac{1}{2}}K_{1}\delta _{1}^{\frac{1}{2}}L_{1}^{\frac{1}{4}}K_{2}^{\frac{1}{%
2}}\left| \Delta \mathbf{q}\right| \\
&\leq &c_{4}^{2}\widetilde{C}_{0}K_{1}^{2}K_{2}L_{1}^{\frac{3}{2}}\delta
_{1}^{2}\frac{2}{\nu }+\frac{\nu }{8}\left| \Delta \mathbf{q}\right| ^{2},
\end{eqnarray*}
and for the fourth we have 
\begin{equation*}
\left| \left( \mathbf{Qf,\Delta q}\right) \right| \leq \frac{2\left| \mathbf{%
Qf}\right| ^{2}}{\nu }+\frac{\nu }{8}\left| \Delta \mathbf{q}\right| ^{2}.
\end{equation*}
We denote 
\begin{eqnarray*}
\frac{1}{2}C_{1}^{2} &=&\frac{2c_{1}^{2}}{\nu }\rho _{0}\rho
_{2}K_{1}^{2}L_{1}\delta _{1}+\frac{2c_{4}^{2}}{\nu }\widetilde{C}_{0}\rho
_{1}\rho _{2}K_{1}L_{1}^{\frac{1}{2}}\delta _{1}^{\frac{3}{2}}+ \\
&&+\frac{2c_{4}^{2}}{\nu }\widetilde{C}_{0}K_{1}^{2}K_{2}L_{1}^{\frac{3}{2}%
}\delta _{1}^{2}+\frac{2\left| \mathbf{Qf}\right| ^{2}}{\nu }.
\end{eqnarray*}

Then the differential inequality for $\left\| \mathbf{q}\right\| $ becomes

\begin{equation*}
\frac{d\left\| \mathbf{q}\right\| ^{2}}{dt}+\nu \Lambda \left\| \mathbf{q}%
\right\| ^{2}\leq C_{1}^{2},
\end{equation*}
that yields 
\begin{equation*}
\left\| \mathbf{q}\left( t\right) \right\| ^{2}\leq \left\| \mathbf{q}\left(
0\right) \right\| ^{2}e^{-\nu \Lambda t}+\frac{1}{\nu \lambda }%
C_{1}^{2}\delta .
\end{equation*}

Let $t_{4}\left( R\right) $ such that for $t\geq t_{4}\left( R\right) \geq
t_{3}\left( R\right) $ the inequality $\left\| \mathbf{q}\left( 0\right)
\right\| ^{2}e^{-\nu \Lambda t}\leq \frac{1}{\nu \lambda }C_{1}^{2}\delta $
holds. For $t\geq t_{4}\left( R\right) $ (\ref{ref norm1q}) holds with $%
\widetilde{C}_{1}=\sqrt{\frac{2}{\nu \lambda }}C_{1}.$

We remark that $L_{1}\delta _{1}=L\left( n\right) \delta \left( n\right)
,\;\ L_{1}^{\frac{1}{2}}\delta _{1}^{\frac{3}{2}}=L\left( n\right) ^{\frac{1%
}{2}}\delta \left( n\right) ^{\frac{3}{2}},\;$ $L_{1}^{\frac{3}{2}}\delta
_{1}^{2}=$ $L\left( n\right) ^{\frac{3}{2}}\delta \left( n\right) ^{2}$ have
values less than 1 for $n=2,$ decrease when $n$ increases for $n\geq 2$ and
tend to zero when $n\rightarrow \infty .$ Hence, $C_{1}$ may be replaced
with a coefficient that depends only on $\,\nu ,\;\lambda ,\;\left| \mathbf{%
Qf}\right| $ and not on $n.$

Moreover, for $n$ large enough, each of the first four terms of $C_{1}$
becomes smaller than $\frac{\left| \mathbf{Qf}\right| }{\nu \sqrt{\lambda }}%
. $

As for the solution $\mathbf{u}$ in \cite{T-N-S}, it can be proved that $%
\mathbf{q}\left( t\right) $ is analytic in time and is the restriction to
the real axis of an analytic function of complex variable defined on a
neighborhood of the real axis, and by using the Cauchy formula, we obtain (%
\ref{ref q prim}).

Finally, from (\ref{eqqu}) we have 
\begin{equation*}
\Delta \mathbf{q=}\frac{1}{\nu }\left[ \frac{d\mathbf{q}}{dt}+\mathbf{QB}(%
\mathbf{p+q})-\mathbf{Qf}\right]
\end{equation*}
and with the above estimates we obtain (\ref{ref delta q}).$\square $

\section{The new modified Galerkin method}

Let us fix a $T>t_{4}\left( R\right) $. The interval $[0,T]$ is the interval
on which we seek the approximate solution. Obviously, all the above
inequalities are valid for $t\in \left[ t_{4}\left( R\right) ,T\right] .$

In this section we just present the method, while in the following section
we estimate the error of the method.

\subsection{The first level}

The first level of our method is related to the post-processed Galerkin
method of \cite{ANT}.

Let $\ \mathbf{p}_{0}\left( t,\mathbf{x}\right) $ be the solution of the
equation (the Galerkin approximation of (\ref{eqpu})): 
\begin{eqnarray}
\mathbf{p}_{0}^{\prime }-\nu \mathbf{\Delta p}_{0}+\mathbf{PB}\left( \mathbf{%
p}_{0}\right) &=&\mathbf{Pf,}  \label{p0} \\
\mathbf{p}_{0}(0) &=&\mathbf{Pu}_{0},  \notag
\end{eqnarray}
and 
\begin{equation*}
\mathbf{q}_{0}(t)\mathbf{=\Phi }_{0}\left( \mathbf{p}_{0}(t)\right) ,
\end{equation*}
where $\mathbf{\Phi }_{0}:\mathbf{P}\mathcal{H}\rightarrow \mathbf{Q}%
\mathcal{H}\;$is the function whose graph is the a.i.m. $\mathcal{M}_{0}\;$%
defined in \cite{FMT}, that is

\begin{equation}
\mathbf{\Phi }_{0}\left( \mathbf{p}\right) =\left( \nu \mathbf{A}\right)
^{-1}\left[ \mathbf{Qf-QB}\left( \mathbf{p}\right) \right] .  \label{q0}
\end{equation}

We define the corresponding approximate solution for the Navier-Stokes
problem (\ref{u-tilda})-(\ref{ci-u-tilda}) as 
\begin{equation}
\mathbf{u}_{0}\left( t\right) =\mathbf{p}_{0}\left( t\right) +\mathbf{q}%
_{0}\left( t\right) .
\end{equation}

Unlike the method of \cite{ANT}, we compute $\mathbf{q}_{0}$ at every moment
of time, and not only at the end of the time interval, $T\;$(in the course
of the numerical implementation of this method, $\mathbf{q}_{0}$ will be
computed at every point of the grid on $[0,T],$ constructed for the
integration of (\ref{p0})) .

\textbf{Remark. }In\textbf{\ }\cite{T-ind-tr} the function $\mathbf{q}%
_{0,m}\left( t\right) =\mathbf{\Phi }_{0}\left( \mathbf{p}\left( t\right)
\right) $ is defined (with $\mathbf{p}\left( t\right) \;$the $\mathbf{P}$
projection of the exact solution), and then the function $\mathbf{u}%
_{0,m}\left( t\right) =\mathbf{p}\left( t\right) +\mathbf{q}_{0,m}\left(
t\right) ,$ is constructed, it's positive trajectory being named ''an
induced trajectory''. Since $\mathbf{p}_{0}\left( t\right) $ is an
approximation of $\mathbf{p}\left( t\right) $, as is proved in Section 7, it
follows that $\left\{ \mathbf{u}_{0}\left( t\right) ;\;t\geq 0\right\} ,$ is
an approximation of this first induced trajectory of \cite{T-ind-tr}.

\subsection{The second level}

The next level is different from both the nonlinear Galerkin methods and the
post-processed Galerkin method already defined in literature. At this level
we make use of $\mathbf{q}_{0}$ calculated in the previous step and define $%
\ \mathbf{p}_{1}$ as the solution of the equation: 
\begin{eqnarray}
\mathbf{p}_{1}^{\prime }-\nu \mathbf{\Delta p}_{1}+\mathbf{PB}\left( \mathbf{%
p}_{1}\mathbf{+q}_{0}\right) &=&\mathbf{Pf},  \label{def-p1} \\
\mathbf{p}_{1}(0) &=&\mathbf{Pu}_{0}.  \notag
\end{eqnarray}
We expect this correction of $\mathbf{p}_{0}$ to be closer to $\mathbf{p}$
than $\mathbf{p}_{0}$ itself. Then we set 
\begin{eqnarray}
\mathbf{q}_{1}\left( t\right) &=&\left( \nu \mathbf{A}\right) ^{-1}\left[ 
\mathbf{Qf}-\mathbf{QB}\left( \mathbf{p}_{1}\left( t\right) \right) -\mathbf{%
QB}\left( \mathbf{p}_{1}\left( t\right) \mathbf{,q}_{0}\left( t\right)
\right) \right. -  \notag \\
&&-\left. \mathbf{QB}\left( \mathbf{q}_{0}\left( t\right) ,\mathbf{p}%
_{1}\left( t\right) \right) \right] .  \label{q1}
\end{eqnarray}

We define the approximate solution for (\ref{u-tilda})-(\ref{ci-u-tilda}) at
this level by 
\begin{equation}
\mathbf{u}_{1}(t)=\mathbf{p}_{1}(t)+\mathbf{q}_{1}(t).  \label{u1}
\end{equation}

\textbf{Remarks. }

\textbf{1. }In the non-linear Galerkin method \cite{DM}, for the
approximation $\mathbf{p}_{1}$ of $\mathbf{p}$, an equation, similar to (\ref
{def-p1}), but\ with $\mathbf{PB}\left( \mathbf{p}_{1}\mathbf{+\Phi }%
_{0}\left( \mathbf{p}_{1}\right) \right) $ instead of $\mathbf{PB}\left( 
\mathbf{p}_{1}\mathbf{+q}_{0}\right) $\ is\ considered.

\textbf{2. }In what concerns $\mathbf{q}_{1}$, the right hand side of (\ref
{q1}) (let us denote it by $\widetilde{\Phi }_{1}(\mathbf{p}_{1}\left(
t\right) ,\mathbf{q}_{0}\left( t\right) )\;)$ was inspired from the function 
$\mathbf{q}_{1,m}\left( t\right) $ of \cite{T-ind-tr}.

This one is defined as 
\begin{eqnarray*}
\mathbf{q}_{1,m}\left( \mathbf{t}\right) &=&\left( \nu \mathbf{A}\right)
^{-1}\left[ \mathbf{Qf}-\mathbf{QB}\left( \mathbf{p}\left( t\right) \right) -%
\mathbf{QB}\left( \mathbf{p}\left( t\right) \mathbf{,q}_{0,m}\left( t\right)
\right) \right. - \\
&&\left. -\mathbf{QB}\left( \mathbf{q}_{0,m}\left( t\right) ,\mathbf{p}%
\left( t\right) \right) \right]
\end{eqnarray*}
and with it's help the function $\mathbf{u}_{1,m}\left( \mathbf{t}\right) =\;%
\mathbf{p}\left( t\right) +\mathbf{q}_{1,m}\left( t\right) $, is defined,
that generates a new induced trajectory. The construction of the second
a.i.m. in \cite{T-ind-tr}, $\mathcal{M}_{1},\;$is based upon the definition
of \ function $\mathbf{q}_{1,m}\left( \mathbf{t}\right) .$ $\mathcal{M}_{1}$
is the graph of a function $\mathbf{\Phi }_{1}:\mathbf{P}\mathcal{H}%
\rightarrow \mathbf{Q}\mathcal{H}$, given by 
\begin{eqnarray*}
\mathbf{\Phi }_{1}\left( \mathbf{X}\right) &=&\left( \nu \mathbf{A}\right)
^{-1}\left[ \mathbf{Qf}-\mathbf{QB}\left( \mathbf{X}\right) -\mathbf{QB}%
\left( \mathbf{X,\Phi }_{0}\left( \mathbf{X}\right) \right) \right. - \\
&&\left. -\mathbf{QB}\left( \mathbf{\Phi }_{0}\left( \mathbf{X}\right) ,%
\mathbf{X}\right) \right] .
\end{eqnarray*}

So, our function $\mathbf{u}_{1}(t)$ is related to an induced trajectory
and, since this one is related to $\mathcal{M}_{1},$ it is also related to
this a.i.m.

\textbf{3. }In the course of the numerical implementation of the method, $%
\mathbf{q}_{1}$ will be computed in the points of the grid on $[0,T]$, since
it's values in these points will be used at the next level.

\subsection{Inductive definition of the high-order approximations}

Let us consider a $k\in \mathbb{N}.\;$We assume that for every $0\leq j\leq
k+1,$ we already constructed $\mathbf{p}_{j}$ and $\mathbf{q}_{j}$. Now, we
define$\;\mathbf{p}_{k+2}$ as the solution of the problem

\begin{eqnarray}
\mathbf{p}_{k+2}^{\prime }-\nu \mathbf{\Delta p}_{k+2}+\mathbf{PB}\left( 
\mathbf{p}_{k+2}\mathbf{+q}_{k+1}\right) &=&\mathbf{Pf},  \label{p_(k+2)} \\
\mathbf{p}_{k+2}(0) &=&\mathbf{Pu}_{0},  \notag
\end{eqnarray}
with $\mathbf{q}_{k+1}$ defined at the preceding step, and then set $\mathbf{%
q}_{k+2}$ as

\begin{eqnarray}
\mathbf{q}_{k+2} &=&\left( \nu \mathbf{A}\right) ^{-1}\left[ \mathbf{Qf-QB(p}%
_{k+2})-\mathbf{QB(p}_{k+2},\mathbf{q}_{k+1})-\right.  \notag \\
&&\left. -\mathbf{QB(\mathbf{q}}_{k+1}\mathbf{,p}_{k+2})-\mathbf{QB(q}_{k},%
\mathbf{q}_{k})-\mathbf{q}_{k}^{\prime }\right] .  \label{q_k+2}
\end{eqnarray}

Naturally, the corresponding approximate solution of (\ref{u-tilda})-(\ref
{ci-u-tilda}) is defined by 
\begin{equation*}
\mathbf{u}_{k+2}\left( t\right) =\mathbf{p}_{k+2}\left( t\right) +\mathbf{q}%
_{k+2}\left( t\right) .
\end{equation*}

\textbf{Remarks. }

\textbf{1.} The right hand side of (\ref{q_k+2}), that we denote by $%
\widetilde{\Phi }_{k+2}(\mathbf{p}_{k+2}\left( t\right) ,\mathbf{q}%
_{k+1}\left( t\right) ,\mathbf{q}_{k}\left( t\right) ),$ is inspired from
inductive the definition of the function $\mathbf{q}_{k+2,m}\left( t\right) $
of \cite{T-ind-tr}, that is

\begin{eqnarray*}
\mathbf{q}_{k+2,m}\left( t\right) &=&\left( \nu \mathbf{A}\right) ^{-1}\left[
\mathbf{Qf-QB}\left( \mathbf{p}\left( t\right) \right) -\right. \\
&&-\mathbf{QB}\left( \mathbf{p}\left( t\right) ,\mathbf{q}_{k+1,m}\left(
t\right) \right) -\mathbf{QB}\left( \mathbf{q}_{k+1,m}\left( t\right) ,%
\mathbf{p}\left( t\right) \right) - \\
&&\left. -\mathbf{QB}(\mathbf{q}_{k,m}\left( t\right) )-\mathbf{q}%
_{k,m}^{\prime }\left( t\right) \right] .
\end{eqnarray*}
Our functions $\mathbf{u}_{k+2},\;k\geq 0$ are, in fact, approximations of
the functions\ $\mathbf{u}_{k+2,m}\,=\mathbf{p}+\mathbf{q}_{k+2,m}\;$that
generate the induced trajectories in \cite{T-ind-tr}. Our construction
by-passes the construction of a.i.m.s. and is based directly upon that of
the induced trajectories. We can call the sets $\left\{ \mathbf{u}%
_{k+2}(t);\;t\geq 0\right\} $ - approximate induced trajectories.

\textbf{2. }The construction of the high accuracy a.i.m., $\mathcal{M}%
_{k+2}, $ is based upon the definition of the function $\mathbf{q}_{k+2,m}$
of \cite{T-ind-tr}. $\mathcal{M}_{k+2}$ is the graph of $\mathbf{\Phi }%
_{k+2}:\mathbf{P}\mathcal{H\rightarrow }\mathbf{Q}\mathcal{H},$

\begin{center}
\begin{eqnarray}
\mathbf{\Phi }_{k+2}\left( \mathbf{X}\right) &=&\left( \nu \mathbf{A}\right)
^{-1}\left[ \mathbf{Qf-QB}\left( \mathbf{X}\right) -\right.  \label{aim} \\
&&-\mathbf{QB}\left( \mathbf{X},\mathbf{\Phi }_{k+1}\left( \mathbf{X}\right)
\right) -\mathbf{QB}\left( \mathbf{\Phi }_{k+1}\left( \mathbf{X}\right) ,%
\mathbf{X}\right) -  \notag \\
&&\left. -\mathbf{QB}(\mathbf{\Phi }_{k}\left( \mathbf{X}\right) )-\mathbf{%
D\Phi }_{k}\left( \mathbf{X}\right) \Gamma _{k}\left( \mathbf{X}\right) %
\right]  \notag
\end{eqnarray}
\end{center}

\noindent where $\mathbf{D\Phi }_{k}\left( \mathbf{X}\right) $ is the
differential of $\mathbf{\Phi }_{k}\left( \mathbf{X}\right) $, and

\begin{center}
$\Gamma _{k}\left( \mathbf{X}\right) =\nu \mathbf{\Delta X}-\mathbf{PB}%
\left( \mathbf{X+\Phi }_{k}\left( \mathbf{X}\right) \right) +\mathbf{Pf.}$
\end{center}

These a.i.m.s or some variant of these are used in the nonlinear Galerkin
methods, and in the postprocessed high-order nonlinear Galerkin methods.

\textbf{3. }If $k+2$ is the last level we construct, than we may compute $%
\mathbf{q}_{k+2}$ only at the moment of interest ($T$ for example) as in the
postprocessed method of \cite{ANT}.

\section{Estimates of the error of the approximate solutions}

In the proof of the main result of this section, we need the following
result that is a direct consequence of Lemma 1 from \cite{ANT}. We denote by 
$\widehat{\mathbf{v}}_{j,l}^{i}$ the coordinate of the function $\mathbf{v}$
with respect to the eigenfunction $\mathbf{w}_{j,l}^{i}$

\textbf{Lemma }\textit{Let }$\mathbf{G}(s)=\sum\limits_{j,l}\left(
\sum\limits_{i=1}^{4}\widehat{G}_{j,l}^{i}\left( s\right) \mathbf{w}%
_{j,l}^{i}\right) $ \textit{and suppose that } 
\begin{equation*}
\left| \widehat{\mathbf{G}}_{j,l}^{i}\left( s\right) \right| \leq
c_{j,l}^{i},\;\;for\,\;\;0\leq j,l\leq m,\;1\leq i\leq 4.
\end{equation*}

\textit{Then } 
\begin{equation}
\left| \int_{0}^{t}e^{-\nu \left( t-s\right) \mathbf{A}}\mathbf{PG}%
(s)ds\right| \leq \frac{1}{\nu }\left[ \sum\limits_{j,k\leq
m}\sum\limits_{i=1}^{4}\frac{\left( c_{j,l}^{i}\right) ^{2}}{\lambda
_{j,l}^{2}}\right] ^{\frac{1}{2}}.  \label{maj-exp}
\end{equation}

Now we can state and prove our main result.

\textbf{Theorem 2 }\textit{The functions }$\mathbf{u}_{k}\left( t\right)
,\;k\geq 0,$ \textit{defined in the previous section,} \ \textit{represent
approximate solutions of the problem }(\ref{eq_u})-(\ref{c-i}), \textit{and
their accuracy increases with }$k.$ \textit{More precisely, the inequality :}

\begin{equation}
\left| \left( \mathbf{u}-\mathbf{u}_{k}\right) \left( t\right) \right| \leq
C\delta ^{5/4+k/2},\;\;\;\;  \label{err uk}
\end{equation}
\textit{holds \ for every} \ $k\geq 0$ \textit{and for} $t\geq t_{4}\left(
R\right) .$

\bigskip \textbf{Proof }We will prove our assertion by induction.

\textbf{1.} We start with $k=0.\;$In \cite{ANT} the following estimate is
proved, for \newline
$\mathbf{f}\in \left[ L_{per}^{2}(\Omega )\right] ^{2}:$%
\begin{equation}
\left| \left( \mathbf{p-p}_{0}\right) \left( t\right) \right| \leq \frac{C}{%
\Lambda ^{5/4}}=C^{\prime }\delta ^{5/4},  \label{p-p0}
\end{equation}
where $C^{\prime }$ is a constant, large for $\nu $ small. Actually, as can
be seen from \cite{ANT} this $C^{\prime }$ is of the order of the product $%
\widetilde{C}_{0}\widetilde{C}_{1}$, with $\widetilde{C}_{0}$, $\widetilde{C}%
_{1}$ the constants of our Theorem 1. Hence we can assume that $C^{\prime }$
is of the form $K\frac{\left| \mathbf{Qf}\right| ^{2}}{\nu ^{2}\lambda ^{3/2}%
}$ for $n$ great enough, with $K$ a number depending on $T$ \ but not on the
data of the problem (see the proof of our Theorem 1).

Let us observe that $\left| \mathbf{p}_{0}\left( t\right) \right| $ is
bounded for large times. Indeed, 
\begin{eqnarray*}
\left| \mathbf{p}_{0}\left( t\right) \right| &=&\left| \mathbf{p}\left(
t\right) +\mathbf{p}_{0}\left( t\right) -\mathbf{p}\left( t\right) \right|
\leq \left| \mathbf{p}\left( t\right) \right| +\left| \mathbf{p}_{0}\left(
t\right) -\mathbf{p}\left( t\right) \right| \\
&\leq &\rho _{0}+C^{\prime }\delta ^{5/4}=\eta _{0},\text{ for }t\text{
large enough.}
\end{eqnarray*}

The same observation is true for $\left\| \mathbf{p}_{0}\left( t\right)
\right\| $ and for $\left| \Delta \mathbf{p}_{0}\left( t\right) \right| $ 
\begin{eqnarray*}
\left\| \mathbf{p}_{0}\left( t\right) \right\| &\leq &\rho _{1}+C^{\prime
}\delta ^{3/4}=\eta _{1},\; \\
\left| \Delta \mathbf{p}_{0}\left( t\right) \right| &\leq &\rho
_{2}+C^{\prime }\delta ^{1/4}=\eta _{2},\;\text{for }t\text{ large enough.}
\end{eqnarray*}

In order to estimate the various norms of $\left( \mathbf{q}-\mathbf{q}%
_{0}\right) \left( t\right) $ we write 
\begin{equation}
\mathbf{q}=\left( \nu \mathbf{A}\right) ^{-1}\left[ \mathbf{Qf}-\mathbf{QB}%
\left( \mathbf{p+q,p+q}\right) +\frac{d\mathbf{q}}{dt}\right] ,
\label{q-expl}
\end{equation}
subtract from this relation the definition relation of $\mathbf{q}_{0},$ (%
\ref{q0}) and apply $\nu \mathbf{\Delta }$ to the obtained equality. In
norm, we have

\begin{eqnarray}
\left| \nu \mathbf{\Delta }\left( \mathbf{q-q}_{0}\right) \right| &=&\left| 
\mathbf{QB}(\mathbf{p+q})-\mathbf{QB}\left( \mathbf{p}_{0}\right) +\frac{d%
\mathbf{q}}{dt}\right|  \notag \\
&\leq &\left| \mathbf{QB}\left( \mathbf{p-p}_{0},\mathbf{p}\right) \right|
+\left| \mathbf{QB}\left( \mathbf{p}_{0},\mathbf{p-p}_{0}\right) \right| + 
\notag \\
&&+\left| \mathbf{QB}(\mathbf{p,q})\right| +\left| \mathbf{QB}(\mathbf{q,p}%
)\right| +\left| \mathbf{QB}\left( \mathbf{q},\mathbf{q}\right) \right|
+\left| \frac{d\mathbf{q}}{dt}\right| .  \label{q-q0}
\end{eqnarray}

For the first term in the right side, with (\ref{B(u,v)2}) we have:

\begin{eqnarray*}
\left| \mathbf{QB}\left( \mathbf{p-p}_{0},\mathbf{p}\right) \right| &\leq
&c_{2}L^{\frac{1}{2}}\left\| \mathbf{p-p}_{0}\right\| \left\| \mathbf{p}%
\right\| \\
&\leq &CL^{\frac{1}{2}}\rho _{1}\delta ^{3/4}\leq C\rho _{1}\delta ^{1/2},
\end{eqnarray*}
where we used once more the inequality $L^{\frac{1}{2}}\delta ^{1/4}\leq 1.$
Here and in the sequel, $C$ denotes a generic constant (not depending on $m$
but depending on $\nu ,\;\mathbf{f,\;}\lambda )$.

The same estimate holds for the second term. With (\ref{B(u,v)1}), the third
term yields:

\begin{eqnarray*}
\left| \mathbf{QB}(\mathbf{p,q})\right| &\leq &c_{1}\left| \mathbf{p}\right|
^{\frac{1}{2}}\left| \mathbf{\Delta p}\right| ^{\frac{1}{2}}\left\| \mathbf{q%
}\right\| \\
&\leq &C\rho _{0}^{1/2}\rho _{2}^{1/2}\delta ^{1/2},
\end{eqnarray*}
and the fourth 
\begin{eqnarray*}
\left| \mathbf{QB}(\mathbf{q,p})\right| &\leq &c_{1}\left| \mathbf{q}\right|
^{\frac{1}{2}}\left| \mathbf{\Delta q}\right| ^{\frac{1}{2}}\left\| \mathbf{p%
}\right\| \\
&\leq &C\rho _{1}\delta ^{1/2}.
\end{eqnarray*}

By using (\ref{ref q prim}) and all the above inequalities in (\ref{q-q0})
we obtain 
\begin{equation}
\left| \nu \Delta \left( \mathbf{q}\left( t\right) \mathbf{-q}_{0}\left(
t\right) \right) \right| \leq C\delta ^{1/2},
\end{equation}
for $t$ great enough. As consequences 
\begin{equation}
\left\| \mathbf{q}\left( t\right) \mathbf{-q}_{0}\left( t\right) \right\|
\leq C\delta ,\;\;\left| \mathbf{q}\left( t\right) \mathbf{-q}_{0}\left(
t\right) \right| \leq C\delta ^{3/2}.  \label{(q-q0)}
\end{equation}

Inequality (\ref{p-p0}) and the second inequality above imply 
\begin{equation}
\left| \mathbf{u}\left( t\right) -\mathbf{u}_{0}\left( t\right) \right| \leq
C\delta ^{5/4}.  \label{err u0}
\end{equation}

We must remark that, as is proved for $\mathbf{u}$ in \cite{T-N-S}, we can
prove that $\mathbf{p}_{0}$ is analytic in time, and more than that, it is
the restriction of an analytic function of a complex variable to the real
axis. This properties are transferred to $\mathbf{q}$ by its definition.
Then, by using the Cauchy formula, it can be proved that 
\begin{equation}
\left| \mathbf{q}^{\prime }\left( t\right) \mathbf{-q}_{0}^{\prime }\left(
t\right) \right| \leq C\delta ^{3/2}.  \label{q-q0prim}
\end{equation}

We also remark, for later use, that (\ref{q0}) and the dissipativity of $%
\mathbf{p}_{0}$ imply \newline
$\left| \Delta \mathbf{q}_{0}\right| \leq C,$ $\left\| \mathbf{q}%
_{0}\right\| \leq C\delta ^{\frac{1}{2}}$ and$\;\left| \mathbf{q}_{0}\right|
\leq C\delta $ for $t\geq t_{2}\left( R\right) .$

\textbf{2. }We now estimate $\left| \mathbf{p}-\mathbf{p}_{1}\right| $ and $%
\left| \mathbf{q}-\mathbf{q}_{1}\right| .$ We have, by subtracting (\ref
{def-p1}) from (\ref{eqpu})

\begin{eqnarray*}
\frac{d}{dt}\left( \mathbf{p-p}_{1}\right) &=&\nu \mathbf{\Delta }\left( 
\mathbf{p}-\mathbf{p}_{1}\right) -\mathbf{PB}\left( \mathbf{p+q-}\left( 
\mathbf{\mathbf{p}_{1}\mathbf{+q}_{0}}\right) \mathbf{,p+q}\right) - \\
&&-\mathbf{PB}\left( \mathbf{p}_{1}\mathbf{+q}_{0},\mathbf{p+q-}\left( 
\mathbf{\mathbf{p}_{1}\mathbf{+q}_{0}}\right) \right) .
\end{eqnarray*}
From here, by using the semigroup of linear operators of infinitesimal
generator $\nu \mathbf{A},$ we obtain 
\begin{eqnarray*}
\frac{d}{dt}e^{\nu t\mathbf{A}}\left( \mathbf{p-p}_{1}\right) \left(
t\right) &=&e^{\nu t\mathbf{A}}\left\{ -\mathbf{PB}\left( \mathbf{p-p}_{1},%
\mathbf{u}\right) -\mathbf{PB}\left( \mathbf{q-q}_{0},\mathbf{u}\right)
-\right. \\
&&\left. -\mathbf{PB}\left( \mathbf{p}_{1}+\mathbf{q}_{0},\mathbf{q-q}%
_{0}\right) -\mathbf{PB}\left( \mathbf{p}_{1}+\mathbf{q}_{0},\mathbf{p-p}%
_{1}\right) \right\} ,
\end{eqnarray*}
and, by integrating 
\begin{eqnarray*}
\left( \mathbf{p-p}_{1}\right) \left( t\right) &=&e^{-\nu t\mathbf{A}}\left( 
\mathbf{p-p}_{1}\right) \left( 0\right) - \\
&&-\int_{0}^{t}e^{-\nu \left( t-s\right) \mathbf{A}}\left\{ \mathbf{PB}%
\left( \mathbf{p-p}_{1},\mathbf{u}\right) +\mathbf{PB}\left( \mathbf{p}_{1}+%
\mathbf{q}_{0},\mathbf{p-p}_{1}\right) \right\} ds- \\
&&-\int_{0}^{t}e^{-\nu \left( t-s\right) \mathbf{A}}\left\{ \mathbf{PB}%
\left( \mathbf{p}_{1}+\mathbf{q}_{0},\mathbf{q-q}_{0}\right) +\mathbf{PB}%
\left( \mathbf{q-q}_{0},\mathbf{u}\right) \right\} ds.
\end{eqnarray*}

Following \cite{ANT} we use the inequalities \cite{CF} 
\begin{equation*}
\left| \mathbf{A}^{-\delta }\mathbf{B}\left( \mathbf{u},\mathbf{v}\right)
\right| \leq \left\{ 
\begin{array}{c}
C\left| \mathbf{A}^{1-\delta }\mathbf{u}\right| \left| \mathbf{v}\right|
\leq C\left| \mathbf{A}^{1/2}\mathbf{u}\right| \left| \mathbf{v}\right| , \\ 
C\left| \mathbf{u}\right| \left| \mathbf{A}^{1-\delta }\mathbf{v}\right|
\leq C\left| \mathbf{u}\right| \left| \mathbf{A}^{1/2}\mathbf{v}\right| ,
\end{array}
\right.
\end{equation*}
valid for $\delta \in \left( 1/2,1\right) $ and \cite{H} 
\begin{equation*}
\left| \mathbf{A}^{\delta }e^{-\nu t\mathbf{A}}\right| \leq Ct^{-\delta }e^{-%
\frac{\nu \lambda }{2}t},
\end{equation*}
and obtain

\begin{center}
$\left| \left( \mathbf{p-p}_{1}\right) \left( t\right) \right| \leq \left|
e^{-\nu t\mathbf{A}}\left( \mathbf{p-p}_{1}\right) \left( 0\right) \right|
+\int_{0}^{t}C\left( t-s\right) ^{-\delta }e^{-\frac{\nu \lambda }{2}\left(
t-s\right) }\left| \left( \mathbf{p-p}_{1}\right) \left( s\right) \right|
ds+ $

$+\left| \int_{0}^{t}e^{-\nu \left( t-s\right) \mathbf{A}}\left[ \mathbf{PB}%
\left( \mathbf{p}_{1}+\mathbf{q}_{0},\mathbf{q-q}_{0}\right) +\mathbf{PB}%
\left( \mathbf{q-q}_{0},\mathbf{p+q}\right) \right] \left( s\right)
ds\right| .$
\end{center}

A form of Gronwall inequality (\cite{H}, Lemma 7.1.1) implies

$\left| \left( \mathbf{p-p}_{1}\right) \left( t\right) \right| \leq $

$\leq C\underset{0\leq t\leq T}{\max }\left| \int_{0}^{t}e^{-\nu \left(
t-s\right) \mathbf{A}}\left\{ \mathbf{PB}\left( \mathbf{p}_{1}+\mathbf{q}%
_{0},\mathbf{q-q}_{0}\right) +\mathbf{PB}\left( \mathbf{q-q}_{0},\mathbf{p+q}%
\right) \right\} \left( s\right) ds\right| .$

We must remark that the constant $C$ above is of the order of $e^{T}.$

By using the method of \cite{ANT}, we find the estimates for the coordinates
of the several terms in the accolade: 
\begin{eqnarray}
\left| \widehat{\mathbf{B}\left( \mathbf{q}_{0},\mathbf{q-q}_{0}\right) }%
_{j,k}\right| &\leq &\left| \mathbf{q}_{0}\right| \left| \mathbf{A}^{\frac{1%
}{2}}\left( \mathbf{q-q}_{0}\right) \right| \leq C\delta \delta =C\delta
^{2},  \label{B(q,q-q0)} \\
\left| \widehat{\mathbf{B}\left( \mathbf{q-q}_{0},\mathbf{q}\right) }%
_{j,k}\right| &\leq &\left| \mathbf{q-q}_{0}\right| \left| \mathbf{A}^{\frac{%
1}{2}}\mathbf{q}\right| \leq C\delta ^{3/2}\delta ^{1/2}=C\delta ^{2},
\label{B(q-q0,q)} \\
\left| \widehat{\mathbf{B}\left( \mathbf{p}_{1},\mathbf{q-q}_{0}\right) }%
_{j,k}\right| &\leq &\left| \mathbf{A}^{\frac{1}{2}}\left( \mathbf{q-q}%
_{0}\right) \right| \left( \left| \left( \mathbf{I}-\mathbf{P}_{m-j}\right) 
\mathbf{p}\right| +\left| \left( \mathbf{I}-\mathbf{P}_{m-k}\right) \mathbf{p%
}\right| \right)  \notag \\
&\leq &C\delta \left( \frac{1}{\lambda _{_{m-j+1}}}+\frac{1}{\lambda
_{_{m-k+1}}}\right) ,  \label{B(p,q-q0)} \\
\left| \widehat{\mathbf{B}\left( \mathbf{q-q}_{0},\mathbf{p}\right) }%
_{j,k}\right| &\leq &CK\delta ^{3/2}\left( \frac{1}{\lambda _{_{m-j+1}}^{%
\frac{1}{2}}}+\frac{1}{\lambda _{_{m-k+1}}^{\frac{1}{2}}}\right) ,  \notag
\end{eqnarray}
where $\mathbf{P}_{m-j}$ represents the projection operator on the space
spanned by the eigenfunctions corresponding to the eigenvalues in $\Gamma
_{m-j}$ and $\lambda _{j}=\lambda _{j,0}.$

By using the inequalities (\ref{maj-exp}), (\ref{B(q,q-q0)}), (\ref
{B(q-q0,q)}) and 
\begin{equation*}
\sum\limits_{j,k\leq m}\lambda _{j,k}^{-2}\leq \widetilde{C},
\end{equation*}
it follows that

\begin{equation*}
\left| \int_{0}^{t}e^{-\nu \left( t-s\right) \mathbf{A}}\left\{ \mathbf{PB}%
\left( \mathbf{q}_{0},\mathbf{q-q}_{0}\right) +\mathbf{PB}\left( \mathbf{q-q}%
_{0},\mathbf{q}\right) \right\} ds\right| \leq \widetilde{C}K\delta ^{2}.
\end{equation*}

In order to estimate the term $\left| \int_{0}^{t}e^{-\nu \left( t-s\right) 
\mathbf{A}}\mathbf{PB}\left( \mathbf{p}_{1},\mathbf{q-q}_{0}\right)
ds\right| $ we use (\ref{B(p,q-q0)}) and the inequality 
\begin{equation*}
\sum\limits_{j,k\leq m}\frac{1}{\lambda _{j,k}^{2}\lambda _{m-j+1}^{2}}\leq 
\frac{C}{\left( m+1\right) ^{3}}=C\delta ^{3/2},
\end{equation*}
proved in \cite{ANT}. It follows

\begin{equation*}
\left| \int_{0}^{t}e^{-\nu \left( t-s\right) \mathbf{A}}\mathbf{PB}\left( 
\mathbf{p}_{1},\mathbf{q-q}_{0}\right) ds\right| \leq C\delta ^{1+\frac{3}{4}%
}.
\end{equation*}

The same estimate can be proved for $\left| \int_{0}^{t}e^{-\nu \left(
t-s\right) \mathbf{A}}\left\{ \mathbf{PB}\left( \mathbf{q-q}_{0},\mathbf{p}%
\right) \right\} ds\right| ,$ hence finally we have 
\begin{equation}
\left| \mathbf{p-p}_{1}\right| \leq C\delta ^{7/4}.  \label{est p-p1}
\end{equation}

We easily see that $\left| \mathbf{p}_{1}\right| \leq \eta _{0},$ $\left\| 
\mathbf{p}_{1}\right\| \leq \eta _{1},\;$ $\left| \Delta \mathbf{p}%
_{1}\right| \leq \eta _{2}.$

Now, in order to estimate the various norms of $\mathbf{q}-\mathbf{q}_{1},$
we subtract (\ref{q1}) from (\ref{q-expl}), we apply the operator $\nu
\Delta $, and take the norm in $\mathcal{H}$ of the resulted equality. After
grouping the terms in a convenient way, we get

\begin{eqnarray}
\left| \nu \mathbf{\Delta }\left( \mathbf{q}-\mathbf{q}_{1}\right) \right|
&\leq &\left| \mathbf{QB}\left( \mathbf{p}-\mathbf{p}_{1},\mathbf{p}\right)
\right| \mathbf{+}\left| \mathbf{QB}\left( \mathbf{p}_{1},\mathbf{p}-\mathbf{%
p}_{1}\right) \right| +  \notag \\
&&+\left| \mathbf{QB}\left( \mathbf{p-p}_{1}\mathbf{,q}\right) \right|
+\left| \mathbf{QB}\left( \mathbf{q}_{0}\mathbf{,p-p}_{1}\right) \right| 
\notag \\
&&+\left| \mathbf{QB}\left( \mathbf{p}_{1}\mathbf{,q-q}_{0}\right) \right|
+\left| \mathbf{QB}\left( \mathbf{q-q_{0},p}\right) \right|  \notag \\
&&+\left| \mathbf{QB}\left( \mathbf{q,q}\right) \right| +\left| \frac{d%
\mathbf{q}}{dt}\right| .  \label{q-q1}
\end{eqnarray}

As we did for $\left| \mathbf{q-q}_{0}\right| $, we estimate one by one the
terms from the right side. For the first one, we use (\ref{B(u,v)2}) and the
inequality $L^{1/2}\delta ^{1/4}\leq 1:$ \ \ 
\begin{eqnarray*}
\left| \mathbf{QB}\left( \mathbf{p}-\mathbf{p}_{1},\mathbf{p,}\right)
\right| &\leq &c_{2}L^{1/2}\left\| \mathbf{p}-\mathbf{p}_{1}\right\| \left\| 
\mathbf{p}\right\| \\
&\leq &CL^{1/2}\delta ^{5/4}\rho _{1}\leq C\delta .
\end{eqnarray*}

The same estimate is valid for the second term. The third term is smaller
than the first and for the fourth the following holds 
\begin{eqnarray*}
\left| \mathbf{QB}\left( \mathbf{q}_{0}\mathbf{,p-p}_{1}\right) \right|
&\leq &c_{1}\left| \mathbf{q_{0}}\right| ^{\frac{1}{2}}\left| \Delta \mathbf{%
q_{0}}\right| ^{\frac{1}{2}}\left\| \mathbf{p-p}_{1}\right\| \\
&\leq &C\delta ^{1/2}\delta ^{5/4}=C\delta ^{7/4}.
\end{eqnarray*}

For the two following terms we use (\ref{B(u,v)2}) respectively (\ref
{b(u,v,w)1}): 
\begin{eqnarray*}
\left| \mathbf{B}\left( \mathbf{p}_{1}\mathbf{,q-q_{0}}\right) \right| &\leq
&c_{1}\left| \mathbf{p}_{1}\right| ^{\frac{1}{2}}\left| \mathbf{\Delta p}%
_{1}\right| ^{\frac{1}{2}}\left\| \mathbf{q-q_{0}}\right\| \\
&\leq &C\eta _{0}^{1/2}\eta _{2}^{1/2}\delta ,
\end{eqnarray*}

\begin{eqnarray*}
\left| \mathbf{B}\left( \mathbf{q-q_{0},p}\right) \right| &\leq &c_{4}\left| 
\mathbf{q-q_{0}}\right| ^{\frac{1}{2}}\left\| \mathbf{q-q_{0}}\right\| ^{%
\frac{1}{2}}\left\| \mathbf{p}\right\| ^{\frac{1}{2}}\left| \Delta \mathbf{p}%
\right| ^{\frac{1}{2}} \\
&\leq &C\delta ^{3/4}\delta ^{1/2}\rho _{1}^{1/2}\rho _{2}^{1/2}\leq C\delta
^{5/4}.
\end{eqnarray*}

The fifth and sixth terms are smaller than the first, respectively the
second term, while for the seventh we have, with (\ref{b(u,v,w)1}) 
\begin{eqnarray*}
\left| \mathbf{B}\left( \mathbf{q,q}\right) \right| &\leq &c_{4}\left| 
\mathbf{q}\right| ^{\frac{1}{2}}\left\| \mathbf{q}\right\| ^{\frac{1}{2}%
}\left\| \mathbf{q}\right\| ^{\frac{1}{2}}\left| \Delta \mathbf{q}\right| ^{%
\frac{1}{2}} \\
&\leq &C\delta .
\end{eqnarray*}

By using the above inequalities in (\ref{q-q1}) we obtain

\begin{equation*}
\left| \nu \mathbf{\Delta }\left( \mathbf{q}-\mathbf{q}_{1}\right) \right|
\leq C\delta .
\end{equation*}
From here 
\begin{equation}
\left\| \mathbf{q}-\mathbf{q}_{1}\right\| \leq C\delta
^{3/2},\;\;\;\;\;\left| \mathbf{q}-\mathbf{q}_{1}\right| \leq C\delta ^{2}.
\label{est q-q1}
\end{equation}

The arguments used to state the analyticity in time of $\mathbf{q_{0}}$
remain valid for $\mathbf{q_{1}}$ and the following relation follows 
\begin{equation*}
\left| \mathbf{q}^{\prime }-\mathbf{q}_{1}^{\prime }\right| \leq C\delta
^{2}.
\end{equation*}
This will be used later. By using (\ref{est p-p1}) and (\ref{est q-q1}) we
now obtain 
\begin{equation}
\left| \mathbf{u}-\mathbf{u}_{1}\right| \leq CL\delta ^{7/4}.  \label{err u1}
\end{equation}
\textbf{3. The induction step.} We assume that, for every $0\leq j\leq k+1$
the inequalities 
\begin{eqnarray*}
\left| \mathbf{p-p}_{j}\right| &\leq &C\delta ^{5/4+j/2}, \\
\left| \mathbf{q-q}_{j}\right| &\leq &C^{\prime }\delta ^{3/2+j/2}, \\
\left| \mathbf{q}^{\prime }\mathbf{-q}_{j}^{\prime }\right| &\leq &C^{\prime
\prime }\delta ^{3/2+j/2} \\
\left\| \mathbf{q-q}_{j}\right\| &\leq &C^{\prime \prime \prime }\delta
^{1+j/2}.
\end{eqnarray*}
hold. We prove that the above inequalities hold also for $j=k+2:$

\begin{center}
$\left( \mathbf{p-p}_{k+2}\right) \left( t\right) =e^{-\nu t\mathbf{A}%
}\left( \mathbf{p-p}_{k+2}\right) \left( 0\right) -$

$-\int_{0}^{t}e^{-\nu \left( t-s\right) \mathbf{A}}\left\{ \mathbf{PB}\left( 
\mathbf{p-p}_{k+2},\mathbf{u}\right) +\mathbf{PB}\left( \mathbf{p}_{k+2}+%
\mathbf{q}_{k+1},\mathbf{p-p}_{k+2}\right) \right\} ds-$

$-\int_{0}^{t}e^{-\nu \left( t-s\right) \mathbf{A}}\left\{ \mathbf{PB}\left( 
\mathbf{p}_{k+2}+\mathbf{q}_{k+1},\mathbf{q-q}_{k+1}\right) +\mathbf{PB}%
\left( \mathbf{q-q}_{k+1},\mathbf{u}\right) \right\} ds.$
\end{center}

As we did for $\left| \left( \mathbf{p-p}_{1}\right) \left( t\right) \right|
,$ we obtain

$\left| \left( \mathbf{p-p}_{k+2}\right) \left( t\right) \right| \leq \left|
e^{-\nu t\mathbf{A}}\left( \mathbf{p-p}_{k+2}\right) \left( 0\right) \right|
+\int_{0}^{t}C\left( t-s\right) ^{-\delta }e^{-\frac{\nu \lambda }{2}\left(
t-s\right) }\left| \mathbf{p-p}_{k+2}\right| ds+$

$+\left| \int_{0}^{t}e^{-\nu \left( t-s\right) \mathbf{A}}\left\{ \mathbf{PB}%
\left( \mathbf{p}_{k+2}+\mathbf{q}_{k+1},\mathbf{q-q}_{k+1}\right) +\mathbf{%
PB}\left( \mathbf{q-q}_{k+1},\mathbf{p+q}\right) \right\} ds\right| .$

The already cited Gronwall-type Lemma of \cite{H} implies

$\left| \left( \mathbf{p-p}_{k+2}\right) \left( t\right) \right| \leq $

$\leq C\underset{0\leq t\leq T}{\max }\left| \int_{0}^{t}e^{-\nu \left(
t-s\right) \mathbf{A}}\left\{ \mathbf{PB}\left( \mathbf{p}_{k+2}+\mathbf{q}%
_{k+1},\mathbf{q-q}_{k+1}\right) +\mathbf{PB}\left( \mathbf{q-q}_{k+1},%
\mathbf{p+q}\right) \right\} ds\right| .$

We evaluate the coordinates of each term in the brackets following $e^{-\nu
\left( t-s\right) \mathbf{A}}$:

\begin{eqnarray*}
\left| \widehat{B\left( \mathbf{q}_{k+1},\mathbf{q-q}_{k+1}\right) }%
_{j,l}\right| &\leq &\left| \mathbf{q}_{k+1}\right| \left| \mathbf{A}^{\frac{%
1}{2}}\left( \mathbf{q-q}_{k+1}\right) \right| \leq C\delta \delta
^{3/2+k/2}=C\delta ^{3/2+\left( k+2\right) /2}, \\
\left| \widehat{B\left( \mathbf{q-q}_{k+1},\mathbf{q}\right) }_{j,l}\right|
&\leq &\left| \mathbf{q-q}_{k+1}\right| \left| \mathbf{A}^{\frac{1}{2}}%
\mathbf{q}\right| \leq C\delta ^{3/2+\left( k+1\right) /2}\delta ^{\frac{1}{2%
}}=C\delta ^{3/2+\left( k+2\right) /2}, \\
\left| \widehat{B\left( \mathbf{p}_{n+2},\mathbf{q-q}_{n+1}\right) }%
_{j,l}\right| &\leq &\left| \mathbf{A}^{\frac{1}{2}}\left( \mathbf{q-q}%
_{n+1}\right) \right| \left( \left| \left( \mathbf{I-P}_{m-j}\right) \mathbf{%
p}_{n+2}\right| +\left| \left( \mathbf{I-P}_{m-l}\right) \mathbf{p}%
_{n+2}\right| \right) \\
&\leq &C\delta ^{3/2+n/2}\left( 1/\lambda _{_{m-j+1}}+1/\lambda
_{_{m-l+1}}\right) , \\
\left| \widehat{B\left( \mathbf{q-q}_{n+1},\mathbf{p}\right) }_{j,l}\right|
&\leq &C\delta ^{3/2+\left( n+1\right) /2}\left( 1/\lambda _{_{m-j+1}}^{%
\frac{1}{2}}+1/\lambda _{_{m-l+1}}^{\frac{1}{2}}\right) .
\end{eqnarray*}

The same arguments used for the terms involved in $\left| \left( \mathbf{p-p}%
_{1}\right) \left( t\right) \right| $ lead to

\begin{eqnarray*}
\left| \int_{0}^{t}e^{-\nu \left( t-s\right) \mathbf{A}}\left\{ \mathbf{PB}%
\left( \mathbf{q}_{k+1},\mathbf{q-q}_{k+1}\right) +\mathbf{PB}\left( \mathbf{%
q-q}_{k+1},\mathbf{q}\right) \right\} ds\right| &\leq &C\delta ^{3/2+\left(
k+2\right) /2}, \\
\left| \int_{0}^{t}e^{-\nu \left( t-s\right) \mathbf{A}}\mathbf{PB}\left( 
\mathbf{p}_{k+2},\mathbf{q-q}_{k+1}\right) ds\right| &\leq &C\delta
^{(3+k)/2+3/4} \\
&=&C\delta ^{(k+2)/2+5/4}.
\end{eqnarray*}

Analogously we can show that 
\begin{equation*}
\left| \int_{0}^{t}e^{-\nu \left( t-s\right) \mathbf{A}}\left\{ \mathbf{PB}%
\left( \mathbf{q-q}_{k+1},\mathbf{p}\right) \right\} ds\right| \leq C\delta
^{(k+2)/2+5/4},
\end{equation*}
and by putting these results together, it follows

\begin{equation}
\left| \left( \mathbf{p-p}_{k+2}\right) \left( t\right) \right| \leq C\delta
^{(k+2)/2+5/4},  \label{p-p(n+2)}
\end{equation}
that confirms our induction hypothesis in what concerns $\mathbf{p}_{k}$. It
also follows that 
\begin{equation*}
\left| \mathbf{p}_{k+2}\right| \leq \eta _{0},\;\;\left\| \mathbf{p}%
_{k+2}\right\| \leq \eta _{1},\;\;\left| \Delta \mathbf{p}_{k+2}\right| \leq
\eta _{2}.
\end{equation*}

Now for $\left| \nu \mathbf{\Delta }\left( \mathbf{q-q}_{k+2}\right) \left(
t\right) \right| $ we have 
\begin{eqnarray*}
\nu \mathbf{\Delta }\left( \mathbf{q}-\mathbf{q}_{k+2}\right) &=&\mathbf{QB}(%
\mathbf{p})\mathbf{-QB(p}_{k+2})+\mathbf{QB}(\mathbf{p,q})\mathbf{-QB(p}%
_{k+2},\mathbf{q}_{k+1})+ \\
&&+\mathbf{QB}(\mathbf{q,p})\mathbf{-QB(q}_{k+1},\mathbf{p}_{k+2})+\mathbf{QB%
}(\mathbf{q,q})\mathbf{-QB(q}_{k},\mathbf{q}_{k})+ \\
&&+\mathbf{q}^{\prime }-\mathbf{q}_{k}^{\prime }
\end{eqnarray*}
and

\begin{eqnarray*}
\left| \nu \mathbf{\Delta }\left( \mathbf{q}-\mathbf{q}_{k+2}\right) \right|
&\leq &\left| \mathbf{QB}(\mathbf{p-p}_{k+2},\mathbf{p})\mathbf{+QB(p}_{k+2},%
\mathbf{p-p}_{k+2})\right| + \\
&&+\left| \mathbf{QB}(\mathbf{p-p}_{k+2}\mathbf{,q}_{k+1})\mathbf{+\mathbf{%
QB(q}_{k+1},\mathbf{p-p}_{k+2})})\right| + \\
&&+\left| \mathbf{QB}(\mathbf{q-q}_{k+1}\mathbf{,p})\mathbf{+\mathbf{QB(p}%
_{k+2},\mathbf{q-q}_{k+1})}\right| + \\
&&+\left| \mathbf{QB}(\mathbf{q-q}_{k}\mathbf{,q})\mathbf{+QB(q}_{k},\mathbf{%
q-q}_{k})\right| +\left| \mathbf{q}^{\prime }-\mathbf{q}_{k}^{\prime
}\right| .
\end{eqnarray*}

We can see, by using the induction hypothesis, (\ref{p-p(n+2)}) and $%
L^{1/2}\delta ^{1/4}\leq 1,\;$that, for the first two terms, we have $\;$ 
\begin{eqnarray*}
\left. 
\begin{array}{c}
\left| \mathbf{QB}(\mathbf{p-p}_{k+2},\mathbf{p})\right| \\ 
\left| \mathbf{QB(p}_{k+2},\mathbf{p-p}_{k+2})\right|
\end{array}
\right\} &\leq &c_{2}L^{1/2}\left\| \mathbf{p-p}_{k+2}\right\| \eta _{1} \\
&\leq &CL^{1/2}\eta _{1}\delta ^{5/4+(k+1)/2} \\
&\leq &C\eta _{1}\delta ^{1/4+(k+2)/2}.
\end{eqnarray*}
The following term is smaller than the first. The fourth term can be
estimated as follows 
\begin{eqnarray*}
\left| \mathbf{\mathbf{QB(q}_{k+1},\mathbf{p-p}_{k+2})}\right| &\leq
&c_{1}\left| \mathbf{q}_{k+1}\right| ^{\frac{1}{2}}\left| \mathbf{\Delta q}%
_{k+1}\right| ^{\frac{1}{2}}\left\| \mathbf{p-p}_{k+2}\right\| \\
&\leq &C\delta ^{1/2}\delta ^{5/4+(k+1)/2}=C\delta ^{5/4+(k+2)/2}.
\end{eqnarray*}

For the fifth term, we obtain 
\begin{eqnarray*}
\left| \mathbf{QB}(\mathbf{q-q}_{k+1}\mathbf{,p})\right| &\leq &c_{4}\left| 
\mathbf{q-q}_{k+1}\right| ^{\frac{1}{2}}\left\| \mathbf{q-q}_{k+1}\right\| ^{%
\frac{1}{2}}\left\| \mathbf{p}\right\| ^{\frac{1}{2}}\left| \Delta \mathbf{p}%
\right| ^{\frac{1}{2}} \\
&\leq &C\rho _{1}^{1/2}\rho _{2}^{1/2}\delta ^{3/4+(k+1)/4}\delta ^{3/4+k/4}
\\
&\leq &C\rho _{1}^{1/2}\rho _{2}^{1/2}\delta ^{3/2+(2k+1)/4},
\end{eqnarray*}
and for the sixth

\begin{eqnarray*}
\left| \mathbf{QB(p}_{k+2},\mathbf{q-q}_{k+1})\right| &\leq &c_{1}\left| 
\mathbf{p}_{k+2}\right| ^{\frac{1}{2}}\left| \mathbf{\Delta p}_{k+2}\right|
^{\frac{1}{2}}\left\| \mathbf{q-q}_{k+1}\right\| \\
&\leq &c_{1}\delta \eta _{0}^{1/2}\eta _{2}^{1/2}\delta ^{3/2+k/2} \\
&\leq &c_{1}\eta _{0}^{1/2}\eta _{2}^{1/2}\delta ^{3/2+(k+2)/2}.
\end{eqnarray*}

Then, by using (\ref{b(u,v,w)2}) we obtain 
\begin{eqnarray*}
\left| \mathbf{QB}(\mathbf{q-q}_{k}\mathbf{,q})\right| &\leq &c_{4}\left| 
\mathbf{q-q}_{k+1}\right| ^{\frac{1}{2}}\left\| \mathbf{q-q}_{k+1}\right\| ^{%
\frac{1}{2}}\left\| \mathbf{q}\right\| ^{\frac{1}{2}}\left| \Delta \mathbf{q}%
\right| ^{\frac{1}{2}} \\
&\leq &C\delta ^{3/4+(k+1)/4}\delta ^{3/4+k/4}\delta ^{1/4}=C\delta
^{1+(k+2)/2}, \\
\left| \mathbf{QB(q}_{k},\mathbf{q-q}_{k})\right| &\leq &c_{4}\left| \mathbf{%
q}_{k}\right| ^{\frac{1}{2}}\left\| \mathbf{q}_{k}\right\| ^{\frac{1}{2}%
}\left\| \mathbf{q-q}_{k}\right\| ^{\frac{1}{2}}\left| \Delta \left( \mathbf{%
q-q}_{k}\right) \right| ^{\frac{1}{2}} \\
&\leq &C\delta ^{1/2}\delta ^{1/4}\delta ^{3/4+k/4}\delta
^{3/4+(k-1)/4}=C\delta ^{1+(k+2)/2}.
\end{eqnarray*}

By using also the induction hypothesis on $\left| \mathbf{q}^{\prime }%
\mathbf{-q}_{k}^{\prime }\right| $ and by comparing the magnitude orders of
the various terms, we find, successively, 
\begin{eqnarray*}
\left| \nu \Delta \left( \mathbf{q}-\mathbf{q}_{k+2}\right) \right| &\leq
&C\delta ^{1/2+(k+2)/2}, \\
\left\| \mathbf{q}-\mathbf{q}_{k+2}\right\| &\leq &C\delta ^{1+(k+2)/2}, \\
\left| \mathbf{q}-\mathbf{q}_{k+2}\right| &\leq &C\delta ^{3/2+(k+2)/2},
\end{eqnarray*}
and these inequalities confirm our induction hypothesis.

From (\ref{p-p(n+2)}) and the above estimates it follows

\begin{equation}
\left| \mathbf{u-u}_{k+2}\right| \leq C\delta ^{5/4+(k+2)/2}.\square
\label{error-k}
\end{equation}

\section{Comments on the method}

\textbf{1. }A major advantage of our method is that we can use very low
dimensional projection spaces for the approximations of $\mathbf{p}$, since
the accuracy of the approximate solution may be increased by using several
iteration levels of the method.

For example, if we choose \textbf{\ }$m=6,$ after having passed through five
levels of the method we obtain an approximate solution $\mathbf{u}_{4}(t)=%
\mathbf{p}_{4}(t)+\mathbf{q}_{4}(t)$ that bears an error of the order of $\
10^{-5}$ since $\delta ^{13/4}=\frac{1}{7^{13/2}}\simeq 0.0000032.$ At each
level $j,\;0\leq j\leq 4$ we will have to solve a system of $4\times
36+4\times 6=168$ ODEs, on the interval $[0,T]$ in order to find the
coordinates of $\mathbf{p}_{j}(t),$ and to compute the coordinates of $%
\mathbf{q}_{j}(t)$ by using algebraic relations.

Here, as in all nonlinear Galerkin methods, problems appear due to $\mathbf{f%
}$. If this function has a infinity of nonzero coefficients in its Fourier
function, it will generate a infinite number of non-zero coordinates in $%
\mathbf{q}_{j}(t).$ A truncation criterion must be applied and it will
depend on $\mathbf{f.\;}$Thus the number of coordinates of $\mathbf{q}%
_{j}(t) $ to be computed depends on $j$ and on the given function $\mathbf{f}
$.

If we chose $m=10,$ we need only four levels of the method for an error of
the order of $10^{-5}$ (in this case $\delta ^{11/4}=\frac{1}{11^{11/2}}%
\simeq 0.00000187).$ But at each level a number of $\ 4\times 100+4\times
10=440$ ODEs must be solved for the coordinates of $\mathbf{p}_{j}(t).$
Besides these, at each level $j$ the coordinates of $\mathbf{q}_{j}(t)$ must
be computed by algebraic relations resulted from the definitions.

\textbf{2.} The program for the integration the systems of ODEs for $\mathbf{%
p}_{j}(t)$ should have the same structure for all $j$, only the coordinates
of $\mathbf{q}_{j-1}(t)$ remaining to be replaced in the nonlinear term.

\textbf{3. }A comparison with the nonlinear Galerkin methods that use
high-accurate a.i.m.s is necessary.

The nonlinear Galerkin method based on the use of high accuracy a.i.m.s \cite
{DM}, \cite{NTW}, applied to the Navier-Stokes problem and corresponding to
our level $k+2,\;k\geq 0$, consists in solving the finite dimensional problem

\begin{eqnarray}
\frac{d\widetilde{\mathbf{p}}}{dt}-\nu \mathbf{\Delta }\widetilde{\mathbf{p}}%
\mathbf{+PB}(\widetilde{\mathbf{p}}+\mathbf{\Phi }_{k+1}\left( \widetilde{%
\mathbf{p}}\right) ) &=&\mathbf{Pf,}  \label{nlG} \\
\widetilde{\mathbf{p}}\left( 0\right) &=&\mathbf{Pu}\left( 0\right) ,  \notag
\end{eqnarray}
for the approximation $\widetilde{\mathbf{p}}\;$of\ $\mathbf{p=Pu.}$ Here,
as above, $\mathbf{\Phi }_{k+1}:\mathbf{P\mathcal{H}\rightarrow Q}\mathcal{H}
$, is the function defining an a.i.m. of high accuracy. The advantage of
this method towards ours is that the system of equations for $\widetilde{%
\mathbf{p}}\;$ is integrated only once. But the problem with solving (\ref
{nlG}) is that the definition of $\mathbf{\Phi }_{k+1}$ requires those of
all $\mathbf{\Phi }_{j}\;$with $j<k+1$ and is very laborious (see \cite{DM}%
). Programming this must be very difficult. The structure of our method,
with iterative levels, makes the computations easier to program, and each
level represents a certain approximation of the solution, so we can enjoy
partial results.

On another hand, all the computations for $\mathbf{q}_{0}\left( t\right) ,\;%
\mathbf{q}_{1}\left( t\right) ,...,\mathbf{q}_{k+1}\left( t\right) $
(together) in our method seem, at first glance, of the same order of
complexity as those necessary for the evaluation of $\mathbf{\Phi }%
_{k+1}\left( \widetilde{\mathbf{p}}\left( t\right) \right) $ in the course
of the numerical integration of (\ref{nlG}) in \cite{DM}, \cite{NTW}.
However, in our method, a major simplification of the computations appears
since $\mathbf{q}_{k-2}^{\prime }$ (from the definition of $\mathbf{q}_{k})$
may be approximated by the numerical derivative $\left( \mathbf{q}_{k-2}(t)-%
\mathbf{q}_{k-2}(t-h)\right) /h$ (since we have already computed $\mathbf{q}%
_{k-2}(t)$ at every time step). This must be compared with the definitions
of\ $z_{j,m}^{\prime }$ in \cite{DM} or $q_{j}^{1}$ in \cite{NTW}, that
yield difficulties in the numerical integration programming. We must also
remark here that the term $\mathbf{D\Phi }_{k-1}\left( \mathbf{X}\right)
\Gamma _{k-1}\left( \mathbf{X}\right) $ in the definition of $\mathbf{\Phi }%
_{k+1}\left( \mathbf{X}\right) $ (that is (\ref{aim}) with $k+1$ instead of $%
k+2)$ is meant to approximate $\mathbf{q}_{k-1,m}^{\prime }$ from the
definition of $\mathbf{q}_{k+1,m}.$ Hence, conceiving a method that uses
directly the functions $\mathbf{q}_{k+1,m}$ instead of the a.i.m.s that are
defined with the help of these functions is very natural.

\textbf{4. }The memory of the computer is better organized \ in our method,
since at the beginning of the computations for the level $\ j$ we may erase
from the memory the value of $\mathbf{p}_{j-1}(t)$ and keep only those of $%
\mathbf{q}_{j-1}(t),\;\mathbf{q}_{j-2}(t).$

\textbf{5. }In order to have not too many computations we may postprocess
our solution $\mathbf{p}_{k+2}(t)$ (at the last level) only at the end of
the time interval $[0,T]$, by adding $\mathbf{q}_{k+2}(T)=\widetilde{\Phi }%
_{k+2}\left( \mathbf{p}_{k+2}(T),\mathbf{q}_{k+1}(T),\mathbf{q}%
_{k}(T)\right) $, as is done in \cite{NTW} for (\ref{nlG}).

\textbf{Acknowledgements}

Work supported by the Romanian Ministry of Education and Research via the
 grant CEEX-05-D11-25/5.11.2005

\end{document}